\documentclass[bj]{imsart}

\RequirePackage[OT1]{fontenc}
\RequirePackage{amsthm,amsmath,amssymb,natbib}
\RequirePackage[colorlinks,citecolor=blue,urlcolor=blue]{hyperref}
\usepackage{latexsym}
\usepackage{subfigure}
\usepackage{mathtools}
\usepackage[linesnumbered,ruled,vlined]{algorithm2e}
\usepackage{algorithmic}
\usepackage{dsfont}

\newcommand{\mP}{\mathbb{P}}
\newcommand{\mE}{\mathbb{E}}

\newcommand{\convAS}{\overset{a.s}{\longrightarrow}}

\newcommand{\convD}{\overset{d}{\longrightarrow}}

\newtheorem{definition}{Definition}[section]
\newtheorem{theorem}{Theorem}[section]
\newtheorem{lemma}{Lemma}[section]

\newtheorem{corollary}{Corollary}[section]

\DeclarePairedDelimiter\floor{\lfloor}{\rfloor}

\arxiv{arXiv:1904.05741}

\startlocaldefs
\numberwithin{equation}{section}
\theoremstyle{plain}

\theoremstyle{remark}
\newtheorem{remark}{Remark}[section]
\endlocaldefs

\begin{document}

\begin{frontmatter}
\title{Comparing a Large Number of Multivariate Distributions}
\runtitle{Comparing a Large Number of Multivariate Distributions}

\begin{aug}
\author{\fnms{Ilmun} \snm{Kim}\thanksref{a,e1}\ead[label=e1,mark]{ilmunk@stat.cmu.edu}}

\address[a]{Department of Statistics and Data Science, Carnegie Mellon University, Pittsburgh, PA 15213, USA
\printead{e1}}

\runauthor{I. Kim}

\affiliation{Some University and Another University}

\end{aug}

\begin{abstract}
	In this paper, we propose a test for the equality of multiple distributions based on kernel mean embeddings. Our framework provides a flexible way to handle multivariate data by virtue of kernel methods and allows the number of distributions to increase with the sample size. This is in contrast to previous studies that have been mostly restricted to classical univariate settings with a fixed number of distributions. By building on Cram{\'e}r-type moderate deviation for degenerate two-sample $V$-statistics, we derive the limiting null distribution of the test statistic and show that it converges to a Gumbel distribution. The limiting distribution, however, depends on an infinite number of nuisance parameters, which makes it infeasible for use in practice. To address this issue, the proposed test is implemented via the permutation procedure and is shown to be minimax rate optimal against sparse alternatives. During our analysis, an exponential concentration inequality for the permuted test statistic is developed which may be of independent interest.
\end{abstract}

\begin{keyword}
\kwd{Bobkov's inequality}
\kwd{Maximum mean discrepancy}
\kwd{K-sample test}
\kwd{Permutation test}
\end{keyword}

\begin{keyword}[class=MSC]
	\kwd{62E20}
	\kwd{62H15}
\end{keyword}

\end{frontmatter}

\section{Introduction}
Let $P_1,\ldots,P_K$ be probability distributions defined on a common measurable space $(\mathcal{X},\mathcal{B})$ for $K \geq 2$. The $K$-sample problem is concerned with testing the null hypothesis $H_0: P_1 = \dots = P_K$ against the alternative hypothesis $H_1: P_i \neq P_j$ for some $i,j \in \{1,\ldots,K\}$. 
This fundamental problem of comparing multiple distributions is a classical topic in statistics with a wide range of applications \citep[][for  reviews]{thas2010comparing,chen2018modern}. Despite its long history, previous approaches to the $K$-sample problem have several limitations. First, many methods are limited to dealing with univariate data. For instance, \cite{kiefer1959k} proposes the $K$-sample analogues of the Kolmogorov--Smirnov and Cram{\'e}r-Von Mises tests. \cite{scholz1987k} generalize the Anderson-Darling test \citep{anderson1952asymptotic} to the $K$-sample case. These approaches are based on empirical distribution functions and are not easily extendable to multivariate data. Some other references that are restricted to the univariate $K$-sample problem include \cite{conover1965several,zhang2007k, wylupek2010data,quessy2012cramer,lemeshko2018some}. Second, most research in this area has been carried out under classical asymptotic regimes where the sample size goes to infinity but the number of distributions is fixed \citep[e.g.,][]{burke1979asymptotic,bouzebda2011k,huvskova2008tests,martinez2008k,jiang2015nonparametric,wang2018nonparametric,armando2018k}. Clearly this classical asymptotic analysis is not appropriate for a dataset with large $K$ and it only provides a narrow picture of the behavior of a test. To the best of our knowledge, \cite{zhan2014testing} is the only study in the literature that considers large $K$. However, their analysis is limited to univariate data with fixed sample size. 
Third, recent developments on the multivariate $K$-sample problem are largely built upon an average difference between distributions \citep{bouzebda2011k,huvskova2008tests,rizzo2010disco,zhan2014testing,wang2018nonparametric,armando2018k}. It is well-known that the test based on an average-type test statistic tends to be powerful against dense alternatives in which many of $P_1,\ldots,P_K$ are different to each other. On the other hand, it tends to suffer from low power against sparse alternatives where only a few of $P_1,\ldots,P_K$ are different from the others. Recently, sparse alternatives have been motivated by numerous applications such as DNA microarray analysis and anomaly detection where there are a small number of treatments that can actually contribute response variables. These applications have led to recent developments of tests tailored to sparse alternatives in the context of testing a high-dimensional vector \citep{jeng2010optimal,fan2015power,liu2020sign}, two-sample mean or covariance testing \citep{cai2013two,cai2014two, cai2014high}, analysis of variance~\citep{arias2011global,cai2014high} and independence testing~\citep{han2017distribution}. To our knowledge, however, a multivariate $K$-sample test specifically designed for sparse alternatives is not available in the current literature.

 

In this study, we propose a new $K$-sample test that addresses the aforementioned limitations of the previous approaches. More specifically, we introduce a $K$-sample test based on the kernel mean embedding method that has been successfully applied to multivariate hypothesis testing. Our test statistic is defined as the maximum of pairwise maximum mean discrepancies \citep{gretton2007kernel,gretton2012kernel}, which leads to a powerful test against sparse alternatives. Throughout this paper, we investigate statistical properties of the proposed test under the asymptotic regime where both the sample size and the number of distributions tend to infinity. 
Below, we summarize our main findings and contributions.
\begin{itemize}	
	\item \textbf{Limiting null distribution:} By building on \cite{drton2018high}, we develop Cram{\'e}r-type moderate deviation for degenerate two-sample $V$-statistics. Based on this result, we study the limiting distribution of the proposed test statistic when the sample size and the number of distributions increase simultaneously. In particular, we show the test statistic converges to a Gumbel distribution under some appropriate conditions. \\[-0.5em]
	\item \textbf{Concentration inequality under permutations:} We demonstrate the usefulness of Bobkov's inequality~\citep{bobkov2004concentration} in studying a concentration inequality for the permuted test statistic. By applying his result, we derive an exponential concentration inequality for the proposed test statistic under permutations. In contrast to usual Hoeffding or Bernstein-type inequalities, the developed inequality relies solely on completely known and easily computable quantities without any moment assumption. \\[-0.5em]
	\item \textbf{Uniform consistency of the permutation test:} Leveraging the developed concentration inequality for the permuted statistic, we prove the uniform consistency of the permutation test over the class of sparse alternatives. Under some regularity conditions, we also show that the power of the permutation test cannot be improved from a minimax point of view.  \\[-0.5em]
	\item \textbf{Empirical power comparison against sparse alternatives:} A simulation study is conducted to compare the performance of the proposed maximum-type test with the existing average-type tests in the literature. The simulation results show that the proposed test consistently outperforms the average-type tests against sparse alternatives based on isotropic Laplace and Gaussian distributions. 
\end{itemize}

\paragraph{Outline.}
The paper is organized as follows. In Section~\ref{Section: Test Statistic}, we briefly review the maximum mean discrepancy and introduce our test statistic. Section~\ref{Section: Limiting distribution} studies the limiting distribution of the proposed test statistic when the sample size and the number of distributions tend to infinity simultaneously. Section~\ref{Section: Permutation Approach} formally introduces permutation procedures. In Section~\ref{Section: Concentration inequalities under permutations}, we provide an exponential concentration inequality for the proposed test statistic under permutations. Section~\ref{Section: Power Analysis} investigates the power of the proposed test and proves its optimality property against sparse alternatives. In Section~\ref{Section: Simulations}, we demonstrate the finite-sample performance of the proposed approach via simulations. Finally, Section~\ref{Section: Conclusion} concludes the paper and discusses future work. The proofs not presented in the main text can be found in Appendix~\ref{Section Appendix}.

\vskip 2em

\section{Test Statistic} \label{Section: Test Statistic}
We start with a brief overview of the maximum mean discrepancy proposed by \cite{gretton2007kernel,gretton2012kernel}. Let $\mathcal{H}$ be a reproducing kernel Hilbert space (RKHS) on $\mathcal{X}$ with a reproducing kernel $h: \mathcal{X} \times \mathcal{X} \mapsto \mathbb{R}$. For two functions $f,g \in \mathcal{H}$, we write the inner product on $\mathcal{H}$ by $\langle f, g \rangle_{\mathcal{H}}$ and the associated norm by $\Vert f \Vert_{\mathcal{H}}$. Given a probability distribution $P$, the kernel mean embedding of $P$ is given by $\mu_h(P)= \mE_{X\sim P} [ h(X,\cdot) ]$. Using the feature map $\psi: \mathcal{X} \mapsto \mathcal{H}$, which satisfies $h(x,y) = \langle \psi(x), \psi(y)\rangle_{\mathcal{H}}$, the kernel mean embedding can also be written as $\mE_{X\sim P}[\psi(X)]$ \citep[see e.g.][for details]{muandet2016kernel}. We now provide the definition of the maximum mean discrepancy (MMD) associated with kernel $h$. 
\begin{definition}[Maximum mean discrepancy]	
	Given two probability distributions, say $P_1$ and $P_2$, such that $\mE_{X_1 \sim P_1} \| \psi(X_1)\|_{\mathcal{H}} <\infty$ and $\mE_{X_2 \sim P_2} \| \psi(X_2) \|_{\mathcal{H}} <\infty$, the maximum mean discrepancy is defined as the RKHS norm of the difference between $\mu_h(P_1)$ and $\mu_h(P_2)$, i.e. 
	\begin{align*}
	\mathcal{V}_h(P_1,P_2) = \| \mu_h(P_1) - \mu_h(P_2) \|_\mathcal{H}.
	\end{align*}
\end{definition}

It has been shown that when kernel $h$ is characteristic~\citep[see e.g.,][]{fukumizu2008kernel,sriperumbudur2011universality}, the MMD becomes zero \emph{if and only if} $P_1 = P_2$. Some examples of characteristic kernels include Gaussian and Laplace kernels on $\mathcal{X} = \mathbb{R}^d$. This characteristic property allows to have a consistent two-sample test against any fixed alternatives. For general $K$-sample cases, we consider the maximum of pairwise maximum mean discrepancies as our metric, i.e. 
\begin{align*}
\mathcal{V}_{h,\text{max}}(P_1,\ldots,P_K) = \max_{1 \leq k < l \leq K}\| \mu_h(P_k) - \mu_h(P_l) \|_\mathcal{H}.
\end{align*}
Hence as long as $h$ is characteristic, it is clear to see that $\mathcal{V}_{h,\text{max}}(P_1,\ldots,P_K)$ is zero \emph{if and only if} $P_1 = \cdots = P_K$. 


Suppose that we observe identically distributed samples $X_{1,k},\ldots,X_{n_k,k} \sim P_k$ for each $k=1,\ldots,K$ and assume that the samples are mutually independent. We propose our test statistic defined as a plug-in estimator of $\mathcal{V}_{h,\text{max}}$:
\begin{align*}
\widehat{\mathcal{V}}_{h,\text{max}} = \max_{1 \leq k < l \leq K} \bigg\Vert \frac{1}{n_{k}} \sum_{i_1=1}^{n_k} \psi(X_{i_1,k}) - \frac{1}{n_{l}} \sum_{i_2=1}^{n_l} \psi(X_{i_2,l})  \bigg\Vert_\mathcal{H}.
\end{align*}
In practice, the test statistic can be computed in a straightforward manner based on the kernel trick \citep[e.g.~Lemma 6 of][]{gretton2012kernel}:
\begin{align*}
\widehat{\mathcal{V}}_{h,\text{max}} = \max_{1 \leq k < l \leq K} \Bigg\{ \frac{1}{n_k^2}\sum_{i_1,i_2=1}^{n_k} h(X_{i_1,k}, X_{i_2,k}) &  + \frac{1}{n_l^2}\sum_{i_1,i_2=1}^{n_l} h(X_{i_1,l}, X_{i_2,l}) \\[.5em]
&  - \frac{2}{n_k n_l} \sum_{i_1=1}^{n_k} \sum_{i_2=1}^{n_l} h(X_{i_1,k}, X_{i_2,l}) \Bigg\}^{1/2}.
\end{align*}
Throughout this paper, we denote the pooled samples by $\{Z_1,\ldots,Z_N\} = \{X_{1,1},\ldots,X_{n_K,K}\}$ where $N = \sum_{k=1}^K n_k$.

\section{Limiting distribution} \label{Section: Limiting distribution}
Given the test statistic, our next step is to determine a critical value of the test with correct size $\alpha$ and good power properties. A common way of calibrating the critical value is based on the limiting null distribution of the test statistic. In this asymptotic approach, the critical value is set to be the $1-\alpha$ quantile of the limiting null distribution and the null hypothesis is rejected when the test statistic exceeds the critical value. The purpose of this section is to demonstrate the difficulty of implementing this asymptotic-based test in our setting. In particular, we show that $\widehat{\mathcal{V}}_{h,\text{max}}$ converges to a Gumbel distribution with a potentially infinite number of unknown parameters under certain conditions. Unfortunately, it is by no means trivial to consistently estimate these infinite nuisance parameters. Furthermore, it is well-known that a maximum-type statistic converges slowly to its limiting distribution \citep[e.g.][]{hall1991convergence}, which also makes the asymptotic test less attractive in practice. These limitations motivate us to delve into the permutation approach later in Section~\ref{Section: Permutation Approach}--\ref{Section: Power Analysis}.

\subsection{Cram{\'e}r-type moderate deviation}
In order to derive the limiting distribution of the maximum of pairwise MMD statistics, it is important to understand the tail behavior of the two-sample MMD statistic. The main tool to this end is Cram{\'e}r-type moderate deviation for degenerate two-sample $V$-statistics that we will develop in this subsection. Our result largely builds upon Cram{\'e}r-type moderate deviation for degenerate one-sample $U$-statistics recently presented by \cite{drton2018high}. 



Let us start with some notation and assumptions. For notational convenience, we write the MMD statistic between $P_1$ and $P_2$ as
\begin{align*}
\widehat{\mathcal{V}}_{12}^2 = \frac{1}{n_1^2}\sum_{i_1,i_2=1}^{n_1} h(X_{i_1,1}, X_{i_2,1}) + \frac{1}{n_2^2}\sum_{i_1,i_2=1}^{n_2} h(X_{i_1,2}, X_{i_2,2})- \frac{2}{n_1 n_2} \sum_{i_1=1}^{n_1} \sum_{i_2=1}^{n_2} h(X_{i_1,1}, X_{i_2,2}).
\end{align*}
By defining $h^\ast(x_1,x_2;y_1,y_2) := h(x_1,x_2) + h(y_1,y_2) - h(x_1,y_1)/2- h(x_1,y_2)/2 - h(x_2,y_1)/2 - h(x_2,y_2)/2$, the MMD statistic can also be written in the form of a two-sample $V$-statistic
\begin{align} \label{Eq: V-statistic}
\widehat{\mathcal{V}}_{12}^2 = \frac{1}{n_1^2n_2^2} \sum_{i_1,i_2=1}^{n_1} \sum_{j_1,j_2=1}^{n_2} h^\ast(X_{i_1,1},X_{i_2,1}; X_{j_1,2},X_{j_2,2}).
\end{align}
Under the null hypothesis, the considered $V$-statistic is \emph{degenerate} meaning that the conditional expectation of $h^\ast(X_{i_1,1},X_{i_2,1}; X_{j_1,2},X_{j_2,2})$ given any one of $X_{i_1,1}, X_{i_2,1}, X_{j_1,2}, X_{j_2,2}$ has zero variance whenever $i_1 \neq i_2$ and $j_1 \neq j_2$.

Let $X_1,X_2$ be independent random vectors from $P_1$. We then define the centered kernel 
\begin{align*}
\overline{h}(x_1,x_2) := h(x_1,x_2) - \mE[h(x_1,X_2)] - \mE[h(X_1,x_2)] + \mE[h(X_1,X_2)],
\end{align*}
which satisfies $\mE[\overline{h}(X_1,X_2)]=0$ and $\mE[\overline{h}(x_1,X_2)]=0$ almost surely. Under the finite second moment condition of the centered kernel, i.e. $\mE[\{\overline{h}(X_1,X_2)\}^2] < \infty$, we may write 
\begin{align} \label{Eq: centered kernel}
\overline{h}(x_1,x_2) = \sum_{v=1}^\infty \lambda_v \varphi_v(x_1) \varphi_v(x_2),
\end{align}
where $\{\lambda_v\}_{v=1}^\infty$ and $\{ \varphi_v(\cdot) \}_{v=1}^\infty$ are the eigenvalues and eigenfunctions of the integral equation $\mE[\overline{h}(x_1,X_2)\varphi_v(X_2)] = \lambda_v \varphi_v(x_1)$ \citep[e.g. page 80 of][]{lee1990u}. 

To facilitate the analysis, we make the following assumptions regarding the kernel function.

\begin{enumerate}\setlength{\itemindent}{0.07in}
	\item[\textbf{(A1).}] Assume that $\mE[|\overline{h}(X_1,X_1)|] < \infty$. \\[-.5em]
	\item[\textbf{(A2).}] Suppose that $\overline{h}(x_1,x_2)$ admits the decomposition in (\ref{Eq: centered kernel}) with $\lambda_1 \geq \lambda_2 \geq \ldots \geq 0$. For all $u,v \in \mathbb{S}^{T-1} := \{x \in \mathbb{R}^T: \|x\|_2 =1 \}$ where $\| \cdot \|_2$ is Euclidean norm in $\mathbb{R}^T$ and any positive integer $T$, assume that there exists a constant $\eta>0$ independent of $T$ such that
	\begin{align} \label{Eq: Multivariate moment condition}
	\mE \big[ \big| \{ \varphi_{1 \cdots T}(X_1)^\top u \}^2 \{\varphi_{1 \cdots T}(X_1)^\top v\}^{m-2} \big| \big]\leq   \eta^m m^{m/2},
	\end{align}
	where $\varphi_{1 \cdots T}(X_1) := (\varphi_1(X_1),\ldots,\varphi_T(X_1))^\top$ and $m=3,4,\ldots$
\end{enumerate}
It is worth noting that the given conditions are more general than those used in \cite{drton2018high}. Specifically, \cite{drton2018high} assume that the kernel $h$ and its eigenfunctions are uniformly bounded. Clearly, \textbf{(A1)} and \textbf{(A2)} are fulfilled under their boundedness assumptions. We also note that $\overline{h}(x_1,x_2)$ is a valid positive definite kernel \citep{sejdinovic2013equivalence}, which yields $\{\overline{h}(x_1,x_2)\}^2 \leq \overline{h}(x_1,x_1) \overline{h}(x_2,x_2)$. Hence, the second moment condition $\mE[\{\overline{h}(X_1,X_2)\}^2] < \infty$ is also satisfied under \textbf{(A1)}. Finally, the multivariate moment condition (\ref{Eq: Multivariate moment condition}) implies that individual eigenfunctions are sub-Gaussian \citep[e.g.][]{vershynin2018high}. 

Under the given conditions, we present Cram{\'e}r-type moderate deviation for the two-sample degenerate $V$-statistic described in (\ref{Eq: V-statistic}). The proof of the following theorem can be found in Appendix~\ref{Section Appendix}.
\begin{theorem}[Cram{\'e}r-type moderate deviation] \label{Theorem: moderate deviation}
	Suppose that \textbf{(A1)} and \textbf{(A2)} are fulfilled. Assume that there exists a constant $C_1 \geq 1$ such that $C_1^{-1} \leq n_1 / n_2 \leq C_1$ and $n_1/ N$ converges to a constant as $N:=n_1+n_2 \rightarrow \infty$. Then under the null hypothesis $P_1=P_2$, we have
	\begin{align} \label{Eq: Cramer-type moderate deviation}
	\frac{\mP ( n_1n_2\widehat{\mathcal{V}}_{12}^2 / N \geq x )}{\mP \left(  \sum_{v=1}^\infty \lambda_v  \xi_v^2 \geq x\right)} = 1 + o(1),
	\end{align}
	uniformly over $x \in (0, o(N^{\theta}))$ where $\xi_1,\xi_2,\ldots$ are independent and identically distributed as $N(0,1)$. Here $\theta$ is a constant that satisfies 
	\begin{align*}
	\theta < \sup \big\{ q \in [0,1/3): \sum\nolimits_{v > \floor{N^{(1-3q)/5}}} \lambda_v = O(N^{-q}) \big\},
	\end{align*}
	when there exist infinitely many non-zero eigenvalues and $\theta = 1/3$ otherwise. 
\end{theorem}


\begin{remark}
	Although we restrict our attention to the two-sample $V$-statistic with a second-order kernel $h^\ast(x_1,x_2;y_1,y_2)$, our result can be straightforwardly extended to higher-order kernels $h^\ast(x_1,\ldots,x_r;y_1,\ldots,y_r)$ for some $r \geq 3$. The key idea is to consider Hoeffding's decomposition of two-sample $U$-statistics \citep[page 40 of][]{lee1990u} and properly control the remainder terms \citep[see,][for one-sample case]{drton2018high}. Finally, using the relationship between $U$- and $V$-statistics \citep[e.g. page 183 of][]{lee1990u}, one can derive the desired result for the $V$-statistic with a higher-order kernel. We do not pursue this direction here since the second-order kernel is enough for our application. 
\end{remark}

\vskip 1em

\subsection{Gumbel limiting distribution} \label{Section: Gumbel limiting distribution}
With the aid of Theorem~\ref{Theorem: moderate deviation}, we are now ready to describe the limiting distribution of the proposed statistic under large $K$ and large $N$ situations. The main ingredient is Chen--Stein method for Poisson approximations \citep{arratia1989two} that has been successfully applied to approximate the distribution of a maximum-type test statistic to a Gumbel distribution \citep[e.g.][]{han2017distribution,drton2018high}. For sake of completeness, we state Theorem 1 of \cite{arratia1989two}.

\begin{lemma}[Theorem 1 of \cite{arratia1989two}] \label{Lemma: Arratia}
	Let $\mathcal{I}$ be an arbitrary index set and for $i \in \mathcal{I}$, let $Y_{i}$ be a Bernoulli random variable with $p_{i} = \mP(Y_i = 1) > 0$. For each $i \in \mathcal{I}$, consider a subset of $\mathcal{I}$ such that $B_i \subset \mathcal{I}$ with $i \in B_i$. Let us define $W = \sum_{i \in \mathcal{I}} Y_{i}$ and $\lambda = \mE(W) = \sum_{i \in \mathcal{I}} p_{i}$. Let $V$ be a Poisson random variable with mean $\lambda$. Then we have that
	\begin{align*}
	\big| \mP(W=0) - \mP(V=0) \big|\leq \min\{1,\lambda^{-1}\} (b_1 + b_2 + b_3)
	\end{align*}
	where 
	\begin{align*}
	b_1 & := \sum_{i \in \mathcal{I}} \sum_{j \in B_{i}} p_{i} p_{j}, 	\ b_2 := \sum_{i \in \mathcal{I}} \sum_{i \neq j \in B_{i}} \mE (Y_i Y_j) \ \text{and} \\[.5em]
	b_3 & := \sum_{i \in \mathcal{I}} \mE \Big| \mE \Big[ Y_{i} - p_{i} \Big| \sum_{j \in \mathcal{I} - B_i} Y_{j} \Big] \Big|.
	\end{align*}
\end{lemma}

Let us denote the two-sample MMD statistic between $P_k$ and $P_l$ by $\widehat{\mathcal{V}}_{kl}^2$, that is $\widehat{\mathcal{V}}_{kl}^2 = \big\Vert n_{k}^{-1} \sum_{i=1}^{n_k} \psi(X_{i,k}) - n_{l}^{-1} \sum_{j=1}^{n_l} \psi(X_{j,l}) \big\Vert_\mathcal{H}^2$. Assume the sample sizes are the same as $n:=n_1=\ldots=n_K$ for simplicity. Then based on the following key observation
\begin{align*}
\mP \big(n\widehat{ \mathcal{V}}_{h,\text{max}}^2/2 \leq x \big) ~ = ~ \mP  \Big\{ \sum\nolimits_{1 \leq k < l \leq K} \mathds{1} \big( n\widehat{\mathcal{V}}_{kl}^2 /2  > x \big) = 0 \Big\},
\end{align*}
Lemma~\ref{Lemma: Arratia} can be applied in our context with $W = \sum_{1 \leq k < l \leq K} \mathds{1} \big( n \widehat{\mathcal{V}}_{kl}^2 / 2  > x \big)$ and $\lambda =  \sum_{1 \leq k < l \leq K} \mP \big( n \widehat{\mathcal{V}}_{kl}^2  / 2 > x \big)$. Ultimately the proof boils down to showing that $b_1,b_2,b_3$ converge to zero under appropriate conditions. This has been established in Appendix~\ref{Section Appendix} and the result is summarized as follows.





\begin{theorem}[Gumbel limit] \label{Theorem: Gumbel Limit}
	Suppose that \textbf{(A1)} and \textbf{(A2)} are fulfilled. Consider a balanced sample case such that $n:=n_1=\ldots=n_K$. Let $\theta$ be a constant chosen as in Theorem~\ref{Theorem: moderate deviation} and assume that $\log K = o(n^{\theta})$. Then under the null hypothesis $P_1=\ldots =P_K$, for any $y \in \mathbb{R}$, 
	\begin{align*}
	&\lim_{n,K\rightarrow \infty }\mP \left( \frac{n}{2\lambda_1} \widehat{\mathcal{V}}^2_{h, \text{\emph{max}}}  - 4\log K -  (\mu_1 - 2) \log \log K   \leq  y\right) \\[.5em]
	= ~ & \exp\Bigg\{ - \frac{2^{\mu_1/2-2}\kappa}{\Gamma(\mu_1/2)}  \exp\left( -\frac{y}{2}\right) \Bigg\}, 
	\end{align*}
	where $\kappa = \prod_{v=\mu_1+ 1}^\infty(1 - \lambda_v/\lambda_1)^{-1/2}$ and $\mu_1$ is the multiplicity of the largest eigenvalue among the sequence $\{\lambda_v\}_{v=1}^\infty$.
\end{theorem}

\begin{remark}
	From Theorem~\ref{Theorem: Gumbel Limit}, it is clear that we need to know or at least estimate a potentially infinite number of parameters $\{\lambda_v\}_{v=1}^\infty$ in order to implement the asymptotic test. Even if one has access to these eigenvalues, the asymptotic test might suffer from slow convergence. This means that the test can be too liberal or too conservative in finite sample size situations.
\end{remark}

\begin{remark}
	When the sample sizes are unbalanced, the limiting distribution of $\widehat{\mathcal{V}}^2_{h, \text{{max}}}$ may not have an explicit expression as in Theorem~\ref{Theorem: Gumbel Limit}. In particular, it depends on the limit values of $n_k /(n_k + n_l)$ for $1\leq k < l \leq K$. To avoid this complication, we simply focus on the case of equal sample sizes and present the explicit formula for the limiting distribution. Nevertheless, if we instead use the weighted $K$-sample statistic:
	\begin{align*}
	\max_{1 \leq k < l \leq K} \bigg( \frac{n_k n_l}{n_k+n_l}\widehat{\mathcal{V}}_{kl}^2 \bigg),
	\end{align*}
	we may obtain the same Gumbel limiting distribution as in Theorem~\ref{Theorem: Gumbel Limit} for general sample sizes.
\end{remark}


\subsection{Examples}

In general, it is challenging to find closed-form expressions for $\{\lambda_v\}_{v=1}^\infty$ and $\{\varphi_v(\cdot)\}_{v=1}^\infty$ as they depend on the kernel as well as the underlying distribution. We end this section with two simple examples for which $\{\lambda_v\}_{v=1}^\infty$ and $\{\varphi_v(\cdot)\}_{v=1}^\infty$ are explicit. Based on these, we illustrate Theorem~\ref{Theorem: Gumbel Limit}.

\begin{itemize}
	\item \textbf{Linear kernel:} Suppose that $\{X_{1,1},\ldots,X_{n,1}, \ldots X_{1,K},\ldots, X_{n,K} \}$ are independent and identically distributed as a multivariate normal distribution with mean zero and covariance matrix $\Sigma$. Suppose further that $\Sigma$ is a diagonal matrix whose diagonal entries are $\lambda_1 = \ldots = \lambda_{\mu_1} > \lambda_{\mu_1+1} \geq \ldots \geq \lambda_d >0$ for some $\mu_1 \geq 1$. Let us consider the linear kernel given as $h(x_1,x_2) = x_1^\top x_2$. Then it is straightforward to see that the centered kernel in (\ref{Eq: centered kernel}) has the eigenfunction decomposition as
	\begin{align*}
	\overline{h}(x_1,x_2) = \sum_{v=1}^d \lambda_v \varphi_v(x_1) \varphi_v(x_2) = \sum_{v=1}^d \lambda_v \big( x_1^{(v)} / \sqrt{\lambda_v} \big) \big( x_2^{(v)} / \sqrt{\lambda_v} \big)
	\end{align*}
	where $x_1^{(v)}$ is the $v$th component of $x_1$. Under the given setting, $\{\varphi_1(X_{1,1}),\ldots,\varphi_d(X_{1,1})\}$ are independent and identically distributed as $N(0,1)$. It can be shown that the conditions in Theorem~\ref{Theorem: Gumbel Limit} are satisfied with $\theta = 1/3$ under the Gaussian assumption. Thus the resulting test statistic converges to a Gumbel distribution as in Theorem~\ref{Theorem: Gumbel Limit}. \\[-.5em]
	\item \textbf{Chi-square kernel:} Suppose that $\{X_{1,1},\ldots,X_{n,1}, \ldots X_{1,K},\ldots, X_{n,K} \}$ are independent and identically distributed on a discrete domain $\{1,\ldots,m\}$ with fixed $m$. Let $p_v > 0$ be the probability of observing the value $v$ among $\{1,\ldots,m\}$ and consequently $\sum_{v=1}^m p_v = 1$. Consider the chi-square kernel defined as $h(x_1,x_2) = \sum_{v=1}^m p_v^{-1} \mathds{1}(x_1=v)\mathds{1}(x_2=v)$. Let $A$ be a $(m-1) \times (m-1)$ matrix whose $(v_1,v_2)$ entry is $a_{v_1,v_2} = p_{v_1}^{-1} + p_{m}^{-1}$ if $v_1=v_2$ and $a_{v_1,v_2} = p_{m}^{-1}$ otherwise. Let us define the eigenfunction $\varphi_v(x)$ to be the $v$th row of $A^{1/2} \{\mathds{1}(x=1) - p_1, \ldots, \mathds{1}(x=m-1) - p_{m-1}\}^\top$ for $v = 1,\ldots,m-1$. Then, following the calculation in Theorem 14.3.1 of \cite{lehmann2006testing},
	\begin{align*}
	\overline{h}(x_1,x_2)  = \sum_{v=1}^{m-1} \lambda_v \varphi_v(x_1) \varphi_v(x_2) = \sum_{v=1}^m \frac{\{\mathds{1}(x_1=v) - p_v\} \{\mathds{1}(x_2=v) - p_v\}}{p_v},
	\end{align*}	
	where $\lambda_1=\ldots =\lambda_{m-1}=1$ and $\lambda_v =0$ for $v \geq m$ and the eigenfunctions are bounded. Thus the conditions in Theorem~\ref{Theorem: Gumbel Limit} are satisfied with $\theta = 1/3$ and the resulting test statistic converges to a Gumbel distribution.
\end{itemize}

\section{Permutation Approach} \label{Section: Permutation Approach}


So far we have investigated the limiting null distribution of the proposed test statistic and demonstrated the difficulty of implementing the resulting asymptotic test. To address the issue, we take an alternative approach based on permutations that does not require prior knowledge on unknown parameters. The key advantage of the permutation approach is that it yields a valid level $\alpha$ test (or a size $\alpha$ test via randomization) for any finite sample size and for any number of distributions. This attractive property is true for any type of underlying distributions, provided that $\{Z_1,\ldots,Z_N\}$ are exchangeable under $H_0$. In the following, we briefly describe the original and randomized permutation procedures. The randomized procedure has a computational advantage over the original procedure by considering a random subset of all permutations. 
\begin{itemize}
	\item \textbf{Permutation approach}: Let $\mathcal{B}_N$ be the collection of all possible permutations of $\{1,\ldots,N\}$. For $\boldsymbol{b}=(b_1,\ldots,b_N) \in \mathcal{B}_N$, we denote by $\widehat{\mathcal{V}}_{h,\text{max}}^{(\boldsymbol{b})}$ the test statistic computed based on the permuted dataset $\{Z_{b_1},\ldots,Z_{b_N} \}$. We then clearly have $ \widehat{\mathcal{V}}_{h,\text{max}}^{(\boldsymbol{b}_0)} = \widehat{\mathcal{V}}_{h,\text{max}}$ for $\boldsymbol{b}_0 = (1,\ldots,N)$. The permutation $p$-value is calculated by
	\begin{align} \label{Eq: permutation p-value}
	p_{\text{perm}} = \frac{1}{N!}\sum_{\boldsymbol{b} \in \mathcal{B}_N} \mathds{1}\left( \widehat{\mathcal{V}}_{h,\text{max}}^{(\boldsymbol{b})} \geq \widehat{\mathcal{V}}_{h,\text{max}} \right).
	\end{align}
	It is well-known that $\mP(p_{\text{perm}}\leq  t) \leq t$ for any $0 \leq t \leq 1$ under $H_0$ \citep[e.g. Chapter 15 of][]{lehmann2006testing}. Consequently $\mathds{1}(p_{\text{perm}} \leq \alpha)$ is a valid level $\alpha$ test. \\[-.5em]
	\item \textbf{Randomized version}: For large $N$, it would be beneficial to consider a subset of $\mathcal{B}_N$ and compute the approximated permutation $p$-value. Suppose that $\boldsymbol{b}_1^\prime,\ldots,\boldsymbol{b}_M^\prime$ are sampled uniformly from $\mathcal{B}_N$ with replacement. We then define a Monte-Carlo version of the permutation $p$-value by
	\begin{align} \label{Eq: Monte-Carlo p-value}
	p_{\text{MC}} = \frac{1}{M+1} \bigg\{ 1 + \sum_{i=1}^M \mathds{1}\left( \widehat{\mathcal{V}}_{h,\text{max}}^{(\boldsymbol{b}_i^\prime)} \geq \widehat{\mathcal{V}}_{h,\text{max}} \right)\bigg\}.
	\end{align}
	It can be shown that $\mP(p_{\text{MC}} \leq  t) \leq t$ for any $0 \leq t \leq 1$ under $H_0$ \citep[e.g. Chapter 15 of][]{lehmann2006testing}. Hence $\mathds{1}(p_{\text{MC}} \leq \alpha)$ is a valid level $\alpha$ test as well.
\end{itemize}
Having motivated the permutation approach, we next analyze uniform consistency as well as minimax optimality of the resulting permutation test against sparse alternatives in Section~\ref{Section: Power Analysis}, building on concentration inequalities developed in the following section.

\section{Concentration inequalities under permutations} \label{Section: Concentration inequalities under permutations}
This section develops a concentration inequality for the permuted MMD statistic with an exponential tail bound. The result established here is especially useful for studying the type II error (or the power) of the proposed permutation test in Section~\ref{Section: Power Analysis}. Our result can also be valuable in addressing the computational issue of the permutation test. The permutation approach suffers from  high computational cost as the number of all possible permutations increases very quickly with the sample size. As a result, it is common in practice to use Monte-Carlo sampling of random permutations to approximate the $p$-value of a permutation test. However, in some application areas such as genetic where extremely small $p$-values are of interest, the Monte-Carlo approach still requires heavy computations \citep{knijnenburg2009fewer,he2019permutation}. Our concentration inequality has an exponential tail bound with completely known quantities. Based on this, one can find a sharp upper bound for the permutation $p$-value (or the permutation critical value) without any computational cost for permutations. We discuss this direction in more detail in Remark~\ref{Remark: p-value}.

\subsection{Bobkov's inequality} \label{Section: Bobkov's inequality}
Before we state the main result of this section, we introduce Bobkov's inequality~\citep{bobkov2004concentration}, which is the key ingredient of our proof. To state his result, we need to prepare some notation in advance. Consider a discrete cube given by
\begin{align*}
\mathcal{G}_{N,m} = \{ \boldsymbol{w} = (w_1,\ldots,w_N) \in \{0,1\}^N : w_1 + \ldots w_N = m \}.
\end{align*}
Note that for each $\boldsymbol{w} \in \mathcal{G}_{N,m}$, there are exactly $m(N-m)$ neighbors $\{s_{ij} \boldsymbol{w} \}_{i \in I(\boldsymbol{w}), j \in J(\boldsymbol{w})}$ where $I(\boldsymbol{w}) = \{i \leq N: w_i = 1 \}$ and $J(\boldsymbol{w}) = \{ j \leq N: w_j = 0 \}$ such that $(s_{ij}\boldsymbol{w})_r = w_r$ for $r \neq i,j$ and $(s_{ij}\boldsymbol{w})_i = w_j, (s_{ij}\boldsymbol{w})_j = w_i$. Now for a function $f$ defined on $\mathcal{G}_{N,m}$, the Euclidean length of discrete gradient $\nabla f(\boldsymbol{w})$ is given as
\begin{align*}
|\nabla f(\boldsymbol{w})|^2 = \sum_{i \in I(\boldsymbol{w})} \sum_{j \in J(\boldsymbol{w})} |f(\boldsymbol{w}) - f(s_{ij}\boldsymbol{w})|^2.
\end{align*}
For more details, we refer to \cite{bobkov2004concentration}. Then Bobkov's inequality is stated as follows:

\begin{lemma}[Theorem 2.1 of \cite{bobkov2004concentration}] \label{Lemma: Bobkov's inequality}
	For every real-valued function $f$ on $\mathcal{G}_{N,m}$ and $| \nabla f(\boldsymbol{w}) | \leq \Sigma$ for all $\boldsymbol{w}$,
	\begin{align*}
	\mP_{\boldsymbol{w}} [ f(\boldsymbol{w}) - \mE_{\boldsymbol{w}} \{f(\boldsymbol{w})\} \geq t] \leq \exp \{ -(N+2)t^2/(4 \Sigma^2) \},
	\end{align*}
	where $\mP_{\boldsymbol{w}}(\cdot)$ represents a counting probability measure on $\mathcal{G}_{N,m}$ and $\mE_{\boldsymbol{w}}(\cdot)$ is the expectation associated with $\mP_{\boldsymbol{w}}(\cdot)$.
\end{lemma}

\subsection{Two-Sample Case} \label{Section: Two-Sample Case}
We first focus on the two-sample case. When $K=2$, it is clear that the proposed test statistic becomes the $V$-statistic in \cite{gretton2012kernel} and
\begin{align} \label{Eq: Two-Sample Statistic}
\widehat{\mathcal{V}}_{h,\text{max}} = \frac{N}{n_2} \bigg\Vert \frac{1}{n_1} \sum_{i_1=1}^{n_1} \psi(X_{i,1}) - \frac{1}{N} \sum_{j=1}^N \psi(Z_{j})  \bigg\Vert_{\mathcal{H}} =  \frac{N}{n_1 n_2} \bigg\Vert \sum_{i_1=1}^{n_1} \overline{\psi}(Z_i)\bigg\Vert_{\mathcal{H}},
\end{align}
where $\overline{\psi}(Z_{i_1}) = \psi(Z_{i_1}) - \frac{1}{N} \sum_{j=1}^N \psi(Z_j)$. Recall that $\boldsymbol{b}$ is a $N$-dimensional random vector uniformly distributed over $\mathcal{B}_N$ in the permutation procedure. As before in Section~\ref{Section: Permutation Approach}, we denote the test statistic based on the permuted dataset $\{Z_{b_1},\ldots, Z_{b_N} \}$ by
\begin{align*}
\widehat{\mathcal{V}}_{h,\text{max}}^{(\boldsymbol{b})} :=  \frac{N}{n_1 n_2} \bigg\Vert \sum_{i_1=1}^{n_1} \overline{\psi}(Z_{b_{i_1}})\bigg\Vert_{\mathcal{H}}.
\end{align*}
We also denote the probability law under permutations (conditional on $Z_1,\ldots,Z_N$) by $\mP_{\boldsymbol{b}}(\cdot)$ and the expectation associated with $\mP_{\boldsymbol{b}}(\cdot)$ by $\mE_{\boldsymbol{b}}(\cdot)$.

It should be stressed that in the two-sample case, there exists $\boldsymbol{w} \in \mathcal{G}_{N,n_1}$ corresponding to each $\boldsymbol{b} \in \mathcal{B}_N$ such that  
\begin{align*}
\widehat{\mathcal{V}}_{h,\text{max}}^{(\boldsymbol{b})} =\widehat{\mathcal{V}}_{h,\text{max}}^{[\boldsymbol{w}]} :=  \frac{N}{n_1 n_2} \bigg\Vert \sum_{i=1}^{N} w_i \overline{\psi}(Z_i)\bigg\Vert_{\mathcal{H}}.
\end{align*}
More importantly, both $\widehat{\mathcal{V}}_{h,\text{max}}^{(\boldsymbol{b})}$ and $\widehat{\mathcal{V}}_{h,\text{max}}^{[\boldsymbol{w}]}$ have the same probability law when $\boldsymbol{b}$ and $\boldsymbol{w}$ are uniformly distributed over $\mathcal{B}_N$ and $\mathcal{G}_{N,n_1}$, respectively. In other words, we have
\begin{align*}
\mP_{\boldsymbol{b}} \big\{\widehat{\mathcal{V}}_{h,\text{max}}^{(\boldsymbol{b})} - \mE_{\boldsymbol{b}} (\widehat{\mathcal{V}}_{h,\text{max}}^{(\boldsymbol{b})})  \geq t \big\} = \mP_{\boldsymbol{w}} \big\{\widehat{\mathcal{V}}_{h,\text{max}}^{[\boldsymbol{w}]} - \mE_{\boldsymbol{w}} (\widehat{\mathcal{V}}_{h,\text{max}}^{[\boldsymbol{w}]})  \geq t \big\} \quad \text{for all $t \in \mathbb{R}$.}
\end{align*}
This key observation allows us to apply Bobkov's inequality to obtain a concentration inequality for the permuted test statistic in the following theorem. 
\begin{theorem}[Concentration inequality for two-sample statistic] \label{Theorem: Concentration Inequality}
	For $K=2$, let $\mP_{\boldsymbol{b}}$ be the uniform probability measure over permutations conditional on $\{Z_1,\ldots,Z_N\}$. Let us write $\gamma_{1,2} = n_1n_2/(n_1+n_2)^2$. Further denote $\widetilde{h}(Z_i,Z_j) = h(Z_i,Z_i) + h(Z_j,Z_j) - 2h(Z_i,Z_j) \geq 0$ and
	\begin{align} \label{Eq: Variance}
	\widehat{\sigma}^2 = \frac{1}{N(N-1)} \sum_{ i \neq j=1}^N \widetilde{h} (Z_i,Z_j).
	\end{align}
	Then for all $t > 0$, we have
	\begin{align} \label{Eq: two-sample concentration}
	\mP_{\boldsymbol{b}} \left( \widehat{\mathcal{V}}_{h,\text{\emph{max}}}^{(\boldsymbol{b})} \geq t + \sqrt{ \frac{\widehat{\sigma}^2}{2N \gamma_{1,2}}} \right) \leq \exp \left(-\frac{N \gamma_{1,2}^2 t^2}{2 \widehat{\sigma}^2}\right).
	\end{align}
	\begin{proof}
		From the previous discussion, it suffices to investigate a concentration inequality for $f(\boldsymbol{w}) := \widehat{\mathcal{V}}_{h,\text{max}}^{[\boldsymbol{w}]}$, which is uniformly distributed on $\mathcal{G}_{N,n_1}$. Since Bobkov's inequality holds for $f(\boldsymbol{w})$, all we need to do is to find meaningful bounds of the expected value of $f(\boldsymbol{w})$ and the Euclidean length of $\nabla f(\boldsymbol{w})$. 
		We first bound the expected value of $f(\boldsymbol{w})$. Using the feature map representation of kernel $h$, it is straightforward to see that
		\begin{align} \label{Eq: Identities}
		& \sum_{i=1}^N \Vert \overline{\psi}(Z_i) \Vert_{\mathcal{H}}^2  = - \sum_{i \neq j=1}^N \big\langle \overline{\psi}(Z_i), \overline{\psi}(Z_j) \big\rangle_{\mathcal{H}} = \frac{1}{2N} \sum_{i \neq j=1}^N \widetilde{h}(Z_i,Z_j).
		\end{align}
		Then using Jensen's inequality together with the above identities, 
		\begin{align*}
		\mE_{\boldsymbol{w}} \bigg[  \bigg\Vert \sum_{i=1}^{N} w_i \overline{\psi}(Z_i)\bigg\Vert_{\mathcal{H}} \bigg] 	\leq ~ & \sqrt{ \mE_{\boldsymbol{w}} \bigg[ \sum_{i=1}^N w_i^2 \Big\Vert  \overline{\psi}(Z_i) \Big\Vert_{\mathcal{H}}^2 + \sum_{i \neq j=1}^N w_i w_j \big\langle \overline{\psi}(Z_i), \overline{\psi}(Z_i) \big\rangle_{\mathcal{H}} \bigg] } \\[.5em]
		= ~ & \sqrt{\frac{n_1}{N} \sum_{i=1}^N \Big\Vert  \overline{\psi}(Z_i) \Big\Vert_{\mathcal{H}}^2  + \frac{n_1(n_1-1)}{N(N-1)} \sum_{i \neq j=1}^N \big\langle \overline{\psi}(Z_i), \overline{\psi}(Z_j) \big\rangle_{\mathcal{H}}} \\[.5em]
		= ~ & \sqrt{\frac{n_1n_2}{2N^2(N-1)} \sum_{i \neq j=1}^N \widetilde{h}(Z_i,Z_j)}.
		\end{align*}
		By multiplying the scaling factor $N/(n_1n_2)$ on both sides, we have $\mE_{\boldsymbol{w}}[f(\boldsymbol{w})] \leq \sqrt{\widehat{\sigma}^2/(2N\gamma_{1,2})}$. 
		
		Next we bound $|\nabla f(\boldsymbol{w})|$. Recall the definition of $s_{ij}w$ in Section~\ref{Section: Bobkov's inequality}. Using the triangle inequality, we see that 
		\begin{align*}
		\Bigg| \frac{N}{n_1n_2} \bigg\Vert \sum_{l=1}^N w_l \overline{\psi}(Z_l)\bigg\Vert_{\mathcal{H}} - \frac{N}{n_1n_2} \bigg\Vert \sum_{l=1}^N (s_{ij}w)_l \overline{\psi}(Z_l)\bigg\Vert_{\mathcal{H}} \Bigg| \leq \frac{N}{n_1n_2}  \bigg\Vert   \overline{\psi}(Z_i) -  \overline{\psi}(Z_j)  \bigg\Vert_{\mathcal{H}}.
		\end{align*}
		Based on this observation, one can find $\Sigma$, which is independent of $\boldsymbol{w}$, as
		\begin{align*}
		| \nabla f(\boldsymbol{w}) |^2 \leq \Sigma^2 := &\frac{N^2}{n_1^2 n_2^2}\sum_{1 \leq i < j \leq N} \Big\Vert   \overline{\psi}(Z_i) -  \overline{\psi}({Z}_j)  \Big\Vert_{\mathcal{H}}^2 =  \frac{N^2}{2n_1^2n_2^2} \sum_{i \neq j=1}^N \widetilde{h}(Z_i,Z_j), 
		\end{align*}
		where the last equality uses the identities in (\ref{Eq: Identities}). Now apply Bobkov's inequality with the above pieces to obtain the desired result. 
	\end{proof}
\end{theorem}

\vskip 0.5em

\begin{remark} \label{Remark: Two-sample}
	Before we move on, we make several comments on Theorem~\ref{Theorem: Concentration Inequality}.
	\begin{enumerate}
		\item[(a)] The tail of the given concentration inequality relies solely on the variance term of the kernel. This is in sharp contrast to Hoeffding or Bernstein-type inequalities~\citep[e.g.][]{boucheron2013concentration} that usually depend on the (possibly unknown) range of random variables. \\[-.5em]
		\item[(b)] The given concentration inequality requires no assumption on random variables such as boundedness or more generally sub-Gaussianity. Furthermore it only depends on known and easily computable quantities in practice.  \\[-.5em]
		\item[(c)] For $0 < \alpha < 1$, consider a test function $\phi_2 : \{Z_1,\ldots,Z_N\} \mapsto \{0,1\}$ such that 
		\begin{align*}
		\phi_2 = I\Bigg\{ \widehat{\mathcal{V}}_{h,\text{{max}}} \geq  \sqrt{\frac{2\widehat{\sigma}^2}{N \gamma_{1,2}^2} \log \left(  \frac{1}{\alpha} \right)} + \sqrt{\frac{\widehat{\sigma}^2}{2N \gamma_{1,2}}} \Bigg\}.
		\end{align*}
		As a corollary of Theorem~\ref{Theorem: Concentration Inequality}, it can be seen that $\phi_2$ is a valid level $\alpha$ test whenever $\{Z_1,\ldots,Z_N\}$ are exchangeable.   \\[-.5em]
		\item[(d)] We stress that our test statistic is a degenerate two-sample $V$-statistic. Therefore, the previous studies on concentration inequalities for the permuted simple sum \citep[e.g.][]{chatterjee2007stein,adamczak2016circular,albert2018concentration} cannot be applied in our context.
	\end{enumerate}
\end{remark}

\subsection{Numerical Illustrations} \label{Section: Numerical Illustrations}
We illustrate the usefulness of Theorem~\ref{Theorem: Concentration Inequality} via simulations. First of all, we can use Theorem~\ref{Theorem: Concentration Inequality} to compute an upper bound for the original permutation $p$-value. In detail, suppose that $n_1 = n_2$ with $N=n_1 + n_2$ for simplicity. Then it is straightforward to see that the permutation $p$-value is less than or equal to
\begin{align*}
p_{\text{Bobkov}} := 
\begin{cases}
\exp \Big\{ - \frac{N}{32 \widehat{\sigma}^2} \Big( \widehat{\mathcal{V}}_{h,\text{max}} -\sqrt{\frac{2\widehat{\sigma}^2}{N}} \Big)^2 \Big\},  & \text{if $\widehat{\mathcal{V}}_{h,\text{max}} \geq \sqrt{\frac{2\widehat{\sigma}^2}{N}}$} \\
1, & \text{else}.
\end{cases}
\end{align*}
By the nature of the permutation test, $p_{\text{Bobkov}}$ is a valid $p$-value in any finite sample size, in a sense that $\mP(p_{\text{Bobkov}} \leq \alpha) \leq \alpha$ under $H_0$. Another way of obtaining a finite-sample valid $p$-value is to use an \emph{unconditional} concentration inequality. For example, \cite{gretton2012kernel} employ McDiarmid's inequality \citep{mcdiarmid1989method} to have an concentration inequality for the MMD $V$-statistic with a bounded kernel. Based on Theorem 7 of \cite{gretton2012kernel} under the bounded kernel assumption $0 \leq h(x,y) \leq B$, another valid $p$-value can be obtained as
\begin{align*}
p_{\text{McDiarmid}}  := 
\begin{cases}
\exp \Big\{ -\frac{N}{8B} \Big( \widehat{\mathcal{V}}_{h,\text{max}} - \sqrt{\frac{32B}{N}} \Big)^2 \Big\},  & \text{if $\widehat{\mathcal{V}}_{h,\text{max}} \geq \sqrt{\frac{32B}{N}}$} \\
1, & \text{else}.
\end{cases}
\end{align*}
Both approaches provide exponentially decaying $p$-values in sample size but we should emphasize that $p_{\text{Bobkov}}$ does not require any moment conditions on the kernel. Even if the kernel is bounded, $p_{\text{Bobkov}}$ would be preferred to $p_{\text{McDiarmid}}$ when $\widehat{\sigma}^2$ is much smaller than $B$. This point is illustrated under the following set-up.

\paragraph{Set-up.} We consider two kernels: 1) energy distance kernel $h(x,y) = (\|x\|_2 + \|y\|_2 - \|x-y\|_2)/2$ and 2) linear kernel $h(x,y) = x^\top y$. Although these kernels are unbounded in general, they are bounded when the underlying distributions have compact support. For this purpose, we consider two truncated normal distributions with the different location parameters $\mu_1 = 1$ and $\mu_2 = -1$ and the same scale parameter $\sigma^2=1$. We let both distributions have the same support as $[-5,5]$ so that we can calculate the bound $B$ for each kernel. For each sample size $N$ among $\{100,200,\ldots,900,1000\}$, the experiments were repeated 200 times to estimate the expected values of the $p$-values.

\begin{figure}[!t]
	\centering 
	\begin{minipage}[b]{0.45\textwidth}
		\includegraphics[width=\textwidth]{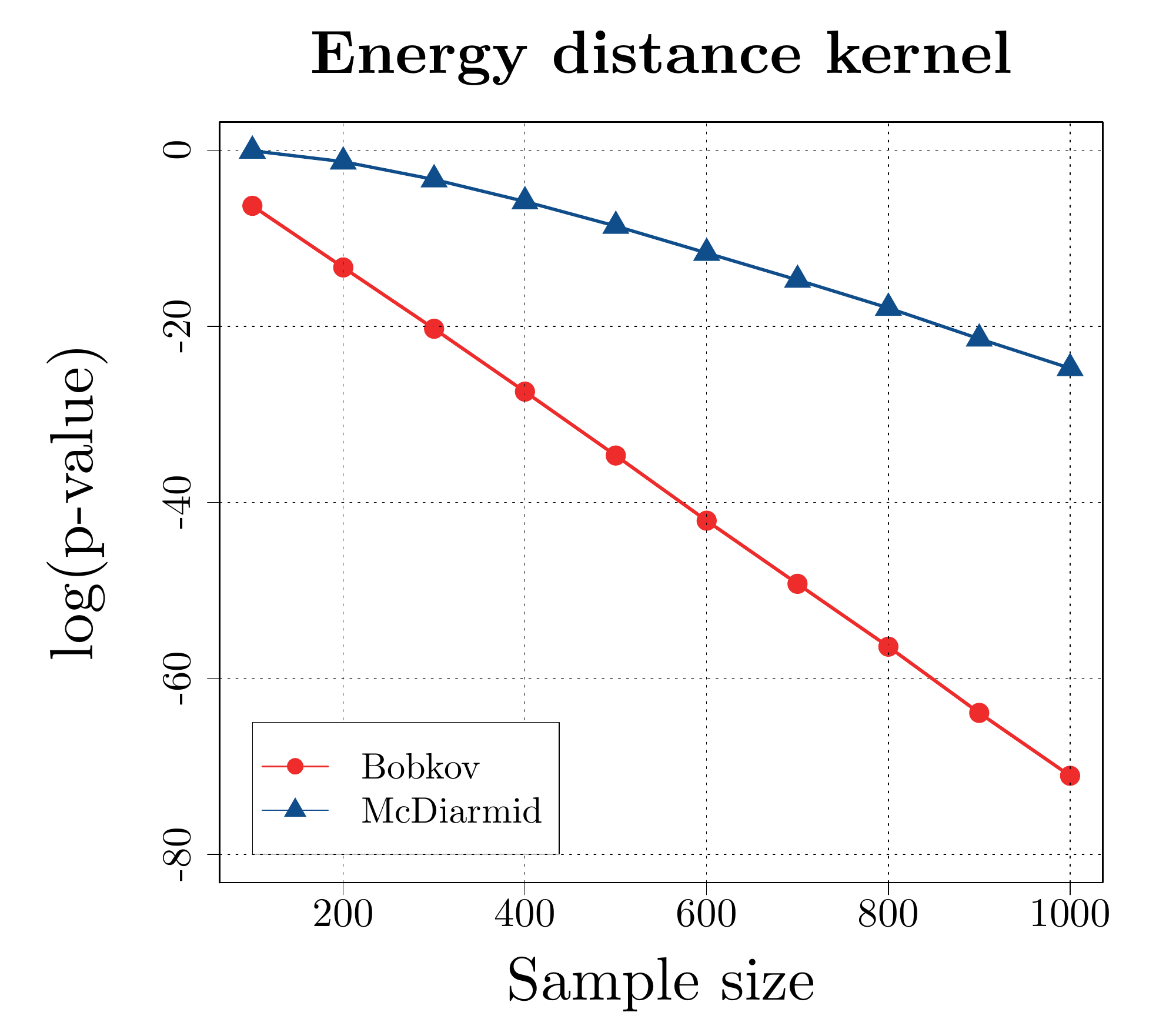}
	\end{minipage}
	\hfill
	\begin{minipage}[b]{0.45\textwidth}
		\includegraphics[width=\textwidth]{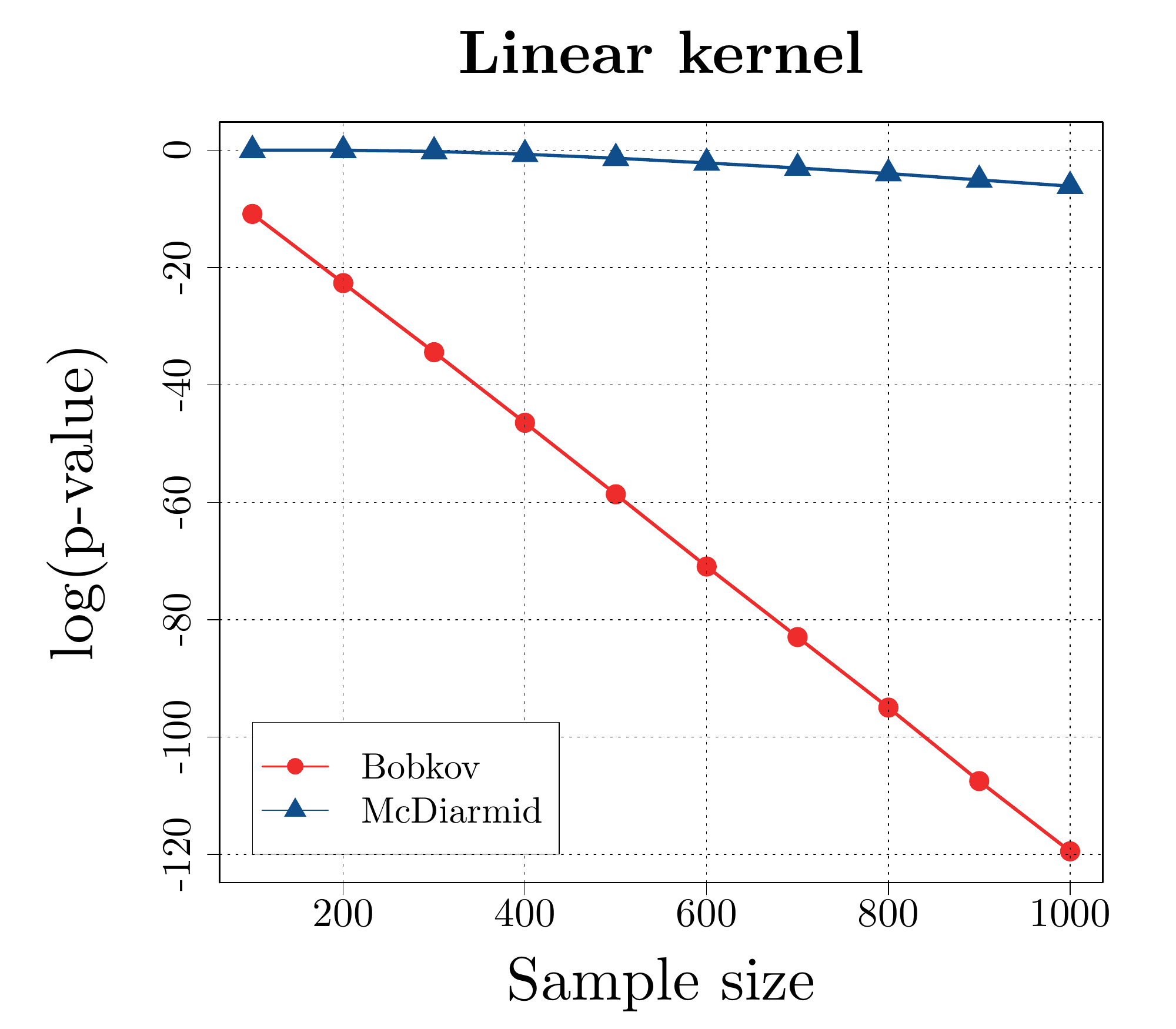}
	\end{minipage} 
	\caption{Comparisons between Bobkov's inequality and McDiarmid inequality in their application to $p$-value evaluation. In both energy distance kernel and linear kernel, Bobkov's inequality returns significantly smaller $p$-values than McDiarmid inequality. See Section~\ref{Section: Numerical Illustrations} for details.} 	\label{Figure: pvalue}
\end{figure}

\paragraph{Results.}
In Figure~\ref{Figure: pvalue}, we present the simulation results of the comparison between $p_{\text{Bobkov}}$ and $p_{\text{McDiarmid}}$ under the described scenario. The $p$-values are displayed in log-scale for better visual comparison. Under the given setting, we observe that $\widehat{\sigma}^2$ is much smaller than $B$ for both kernels, which in turns leads to a smaller value of $p_{\text{Bobkov}}$ compared to $p_{\text{McDiarmid}}$. More specifically, we observe 1) $\widehat{\sigma}^2 \approx 1.61$ on average and $B=10$ for the energy distance kernel and 2) $\widehat{\sigma}^2 \approx 4.01$ on average and $B=100$ for the linear kernel. It is worth noting that the benefit of using $p_{\text{Bobkov}}$ becomes more evident for unbounded random variables for which $p_{\text{McDiarmid}}$ is not even applicable.

\begin{remark} \label{Remark: p-value}
	The test based on $p_{\text{{Bobkov}}}$ may not be recommended when the sample size is small and the significance level $\alpha$ is of moderate size (e.g.~$\alpha = 0.05$). In this case, the permutation test via Monte-Carlo simulations would be more satisfactory. However, when the sample size is large and the significance level is very small (e.g.~$\alpha = 10^{-100})$, the Monte-Carlo approach would be computationally infeasible, requiring at least $\alpha^{-1}$ random permutations in order to reject $H_0$. In this large-sample and small $\alpha$ situation, the approach based on $p_{\text{{Bobkov}}}$ would be practically valuable, which does not require any computational cost on permutations.
\end{remark}

\begin{remark}
	While we focused on the case where $\widehat{\sigma}^2 \ll B$ to highlight the advantage of $p_{\text{Bobkov}}$, it is definitely possible to observe that $p_{\text{McDiarmid}}$ is smaller than $p_{\text{Bobkov}}$, especially when $B$ is comparable to or smaller than $\widehat{\sigma}^2$.
\end{remark}

\allowdisplaybreaks

\subsection{$K$-Sample Case}
Next we give a general result for arbitrary $K \geq 2$. Unfortunately, we cannot directly apply Bobkov's inequality when $K>2$ since the inequality holds only for a function $f(\boldsymbol{w})$ defined on a binary discrete cube. Our strategy to overcome this problem is to first apply Bobkov's inequality to each pairwise MMD test statistic and then aggregate the results via the union bound. To start, we introduce $\widehat{\sigma}_K^2$ in Algorithm~\ref{Alg: sigmaK} that generalizes $\widehat{\sigma}^2$ to the $K$-sample case.

\begin{algorithm}[!htb]
	\caption{Calculation of $\widehat{\sigma}_K^2$} 
	\label{Alg: sigmaK} 
	\begin{algorithmic} \normalsize
		\REQUIRE the pooled samples $\{Z_1,\ldots, Z_N\}$, the number of samples $n_1,\ldots,n_K$.
		\begin{enumerate}\setlength{\itemindent}{-0.15in}
			\item[(1)] Calculate $\widetilde{h}(Z_i,Z_j)$ for $1\leq i \neq j \leq N$. \\[-0.5em]
			\item[(2)] Sort and denote the previous outputs by $\widetilde{h}_{[1]}\geq \ldots \geq \widetilde{h}_{[N(N-1)]}$. \\[-0.5em]
			\item[(3)] Compute $\widehat{\sigma}_K^2 := \max_{1 \leq k < l \leq K} \overline{\sigma}_{kl}^2$ where $\overline{\sigma}_{kl}^2$ is the sample average of $\widetilde{h}_{[1]},\widetilde{h}_{[2]},\ldots$, $\widetilde{h}_{[(n_k+n_l)(n_k+n_l-1)]}$.  \\[-0.5em]
			\item[(4)] Return $\widehat{\sigma}_K^2$.
		\end{enumerate}
	\end{algorithmic}
\end{algorithm}

It can be seen that $\widehat{\sigma}_K^2$ is the same as $\widehat{\sigma}^2$ in (\ref{Eq: Variance}) when $K=2$ and can be computed in quadratic time for large $K$. Using $\widehat{\sigma}_K^2$, we extend Theorem~\ref{Theorem: Concentration Inequality} as follows.


\begin{theorem}[Concentration inequality for $K$-sample statistic] \label{Theorem: K-sample concentration}
	For $K \geq 2$, let $\mP_{\boldsymbol{b}}$ be the uniform probability measure over permutations conditional on $\{Z_1,\ldots,Z_N\}$. For distinct $k,l \in \{1,\ldots,K\}$, let $\gamma_{k,l} = n_k n_l /(n_k+n_l)^2$ and consider $\widehat{\sigma}_K^2$ in Algorithm~\ref{Alg: sigmaK}. Then for any $t \geq 0$,
	\begin{align} \nonumber
	& \mP_{\boldsymbol{b}} \Bigg\{ \widehat{\mathcal{V}}_{h,\text{\emph{max}}}^{(\boldsymbol{b})} \geq t +   \max_{1 \leq k < l \leq K} \sqrt{ \frac{\widehat{\sigma}_K^2}{2(n_k+n_l)\gamma_{k,l}}} \Bigg\} \\[.5em]
	\leq ~ & \binom{K}{2} \exp \Bigg\{ - \min_{1 \leq k < l \leq K} \frac{(n_k+n_l) \gamma_{k,l}^2 t^2}{2\widehat{\sigma}_K^2} \Bigg\}.  \label{Eq: K-sample concentration}
	\end{align}
	\begin{proof}
		For a given permutation $\boldsymbol{b} \in \mathcal{B}_N$, let us denote
		\begin{align*}
		\widehat{\mathcal{V}}_{kl}^{(\boldsymbol{b})} = \bigg\Vert \frac{1}{n_{k}} \sum_{i=1}^{n_k} \psi(Z_{b_{m_{k-1}+i}}) - \frac{1}{n_{l}} \sum_{j=1}^{n_l} \psi(Z_{b_{m_{l-1}+j}})  \bigg\Vert_\mathcal{H},
		\end{align*}
		where $m_{l-1} = \sum_{k=1}^{l-1}n_k$ and $m_0 = 0$ so that $\widehat{\mathcal{V}}_{h,\text{{max}}}^{(\boldsymbol{b})} = \max_{1 \leq k < l \leq K} \widehat{\mathcal{V}}_{kl}^{(\boldsymbol{b})}$. Based on the triangle inequality and the union bound, observe that
		\begin{align}  \nonumber 
		&\mP_{\boldsymbol{b}} \Bigg\{ \widehat{\mathcal{V}}_{h,\text{{max}}}^{(\boldsymbol{b})}  \geq  t  + \max_{1 \leq k < l \leq K} \sqrt{ \frac{\widehat{\sigma}_K^2}{2(n_k+n_l)\gamma_{k,l}}} \Bigg\} \\[.5em] \nonumber 
		\leq ~ & \mP_{\boldsymbol{b}} \left[ \max_{1 \leq k < l \leq K} \Bigg\{ \widehat{\mathcal{V}}_{kl}^{(\boldsymbol{b})} -  \sqrt{ \frac{\widehat{\sigma}_K^2}{2(n_k+n_l) \gamma_{k,l}}} \Bigg\} \geq  t  \right]  \\[.5em]
		\leq ~ & \sum_{1 \leq k < l \leq K} \mP_{\boldsymbol{b}} \Bigg\{  \widehat{\mathcal{V}}_{kl}^{(\boldsymbol{b})} \geq t + \sqrt{ \frac{\widehat{\sigma}_K^2}{2(n_k + n_l) \gamma_{k,l}}} \Bigg\}. \label{Eq: Mid}
		\end{align}
		Let $\widetilde{Z}=\{\widetilde{Z}_1,\ldots,\widetilde{Z}_{n_k+n_l}\}$ be the $n_k+n_l$ samples uniformly drawn from $\{Z_1,\ldots,Z_N\}$ without replacement. Write
		\begin{align*}
		\widehat{\mathcal{V}}_{kl}^{[\boldsymbol{w}]} = \frac{n_k+n_l}{n_k n_l}  \bigg\Vert \sum_{i_1=1}^{n_k+n_l} w_{i_1} \bigg\{ \psi(\widetilde{Z}_{i_1}) - \frac{1}{n_k+n_l} \sum_{i_2=1}^{n_k+n_l} \psi(\widetilde{Z}_{i_2})  \bigg\} \bigg\Vert_{\mathcal{H}},
		\end{align*}
		where $\boldsymbol{w}=\{w_1,\ldots,w_{n_k+n_l}\}$ is a set of Bernoulli random variables uniformly distributed on $ \mathcal{G}_{n_k+n_l,n_k}$ as before. Then by the law of total expectation and a slight modification of the proof of Theorem~\ref{Theorem: Concentration Inequality}, it can be seen that 
		\begin{align*}
		\mP_{\boldsymbol{b}} \Bigg(  \widehat{\mathcal{V}}_{kl}^{(\boldsymbol{b})} \geq t + \sqrt{ \frac{\widehat{\sigma}_K^2}{2(n_k + n_l) \gamma_{k,l}}} \Bigg)  & = ~ \mE_{\widetilde{Z}} \Bigg[ \mP_{\boldsymbol{w}} \Bigg\{ \widehat{\mathcal{V}}_{kl}^{[\boldsymbol{w}]} \geq t + \sqrt{ \frac{\widehat{\sigma}_K^2}{2(n_k + n_l) \gamma_{k,l}}} ~ \Bigg| ~ \widetilde{Z} \Bigg\} \Bigg] \\[.5em]
		& \leq ~ \mE_{\widetilde{Z}} \left[ \exp \Bigg\{-\frac{(n_k+n_l) \gamma_{k,l}^2 t^2}{2\widehat{\sigma}_K^2}\Bigg\} \right] \\[.5em]
		& = ~ \exp \Bigg\{ -\frac{(n_k+n_l) \gamma_{k,l}^2 t^2}{2\widehat{\sigma}_K^2} \Bigg\},
		\end{align*}
		where the last equality follows since $\widehat{\sigma}_K^2$ is invariant to the choice of $\widetilde{Z}$. By putting this result into the right-hand side of (\ref{Eq: Mid}), the proof is complete. 
	\end{proof}
\end{theorem}

\begin{remark}  \label{Remark: K-sample}
	We provide some comments on Theorem~\ref{Theorem: K-sample concentration}.
	\begin{enumerate}
		\item[(a)] When $K=2$, the concentration inequality given in (\ref{Eq: K-sample concentration}) recovers the one in (\ref{Eq: two-sample concentration}). \\[-.5em]
		\item[(b)] One can replace $\widehat{\sigma}_K^2$ with $\max_{1 \leq i < j \leq N} \widetilde{h}(Z_i,Z_j)$ in (\ref{Eq: K-sample concentration}), which takes less time to compute, but at the expense of the loss of the tightness. Note, however, that the bound with $\max_{1 \leq i < j \leq N} \widetilde{h}(Z_i,Z_j)$ is tight enough to prove minimax rate optimality of the proposed test. See the proof of Theorem~\ref{Theorem: Upper Bound} for details. \\[-.5em]
		\item[(c)] As before in the two-sample case, the proposed $K$-sample concentration inequality is valid without any moment condition and it depends solely on known and easily computable quantities.   \\[-.5em]
		\item[(d)] Consider a test function $\phi_K: \{Z_1,\ldots,Z_N\} \mapsto \{0,1\}$ such that 
		\begin{align*}
		\phi_K = I\Bigg[ \widehat{\mathcal{V}}_{h,\text{{max}}} \geq  & \max_{1 \leq k < l \leq K} \sqrt{  \Bigg\{  \frac{2\widehat{\sigma}_K^2}{(n_k+n_l) \gamma_{k,l}^2} \Bigg\} \log \Bigg\{ \frac{\binom{K}{2}}{\alpha} \Bigg\}  } \\[.5em]
		&~~~~~~~~~~~~~~~ +  \max_{1 \leq k < l \leq K} \sqrt{ \frac{ \widehat{\sigma}_K^2 }{2(n_k+n_l)\gamma_{k,l}}} \Bigg].
		\end{align*}
		As a corollary of Theorem~\ref{Theorem: K-sample concentration}, it can be seen that $\phi_K$ is a valid level $\alpha$ test whenever $\{Z_1,\ldots,Z_N\}$ are exchangeable under $H_0$.
	\end{enumerate}
\end{remark}

\section{Power Analysis} \label{Section: Power Analysis}
In this section, we study the power of the permutation test based on the proposed test statistic and prove its minimax rate optimality against certain sparse alternatives. Throughout this section, we need the following assumptions:
\begin{enumerate}\setlength{\itemindent}{0.15in}
	\item[\textbf{(B1).}] Assume that kernel $h$ is uniformly bounded by $0 \leq h(x,y) \leq B$ for all $x,y \in \mathcal{X}$.  \\[-.5em]
	\item[\textbf{(B2).}] There exists a fixed constant $c >0$ such that $n_{\text{max}}/n_{\text{min}} \leq c$ for any sample sizes where $n_{\text{max}}$ and $n_{\text{min}}$ are the maximum and the minimum of $\{n_1,\ldots,n_K\}$ respectively.
\end{enumerate}
Note that the assumption $\textbf{(B1)}$ is satisfied by some widely used kernels e.g. Gaussian and Laplace kernels. It can also be satisfied by many other kernels when the underlying distributions have compact support. 
The second assumption $\textbf{(B2)}$ states that $n_1,\ldots,n_K$ are well-balanced. This assumption, for example, holds for the equal sample sizes with $c=1$. 

\subsection{Power of the permutation test}
Let $\mathcal{P}$ be the set of all distributions on $(\mathcal{X},\mathcal{B})$. We characterize the difference between the null and the alternative in terms of $\max_{1 \leq k < l \leq K} \mathcal{V}_h(P_k,P_l)$, which is the population counterpart of the proposed test statistic $\widehat{\mathcal{V}}_{h,\text{max}}$. In particular, for a given positive sequence $\epsilon_N$ and kernel $h$, let us define a class of alternatives:
\begin{align} \label{Eq: Family of distributions}
\mathcal{F}_h(\epsilon_N) = \big\{(P_1,\ldots,P_K) \in \mathcal{P} : \max_{1 \leq k < l \leq K} \mathcal{V}_{kl} \geq \epsilon_N \big\},
\end{align}
where $\mathcal{V}_{kl} = \mathcal{V}_h(P_k,P_l)$ for simplicity. We call the collection of alternatives in $\mathcal{F}_h(\epsilon_N)$ as the sparse alternatives, in a sense that only a few of $\{\mathcal{V}_{kl}\}_{1 \leq k < l \leq K}$ are required to be greater than $\epsilon_N$ while the rest of them can be zero. Such sparse alternatives have been considered by many authors including \cite{cai2013two}, \cite{cai2014two} and \cite{han2017distribution} in different contexts. The main goal of this subsection is to characterize the conditions under which the permutation test can be uniformly powerful over $\mathcal{F}_h(\epsilon_N)$. More specifically, we show that as long as the lower bound $\epsilon_N$ is sufficiently larger than 
\begin{align*}
r_N^\star := \sqrt{\frac{\log K}{n_{\text{min}}}},
\end{align*}
then the proposed permutation test is uniformly consistent. Furthermore, in Section~\ref{Section: Minimax rate optimality}, we prove that this rate cannot be improved from a minimax perspective under some mild conditions on kernel $h$. In other words, the proposed test is minimax rate optimal against the sparse alternatives with the minimax rate $r_N^\star$.

We start by providing one lemma, which states that $\max_{1 \leq k < l \leq K} |\widehat{\mathcal{V}}_{kl} - \mathcal{V}_{kl}|$ is bounded by $C\sqrt{\log K / n_{\text{min}}}$ for some constant $C$ with high probability.

\begin{lemma} \label{Lemma: Unconditional Analysis}
	Suppose that $\textbf{\emph{(B1)}}$ holds and recall that $\widehat{\mathcal{V}}_{kl} = \big\Vert n_{k}^{-1} \sum_{i_1=1}^{n_k} \psi(X_{i_1,k}) - n_{l}^{-1} \sum_{i_2=1}^{n_l} \psi(X_{i_2,l}) \big\Vert_\mathcal{H}$. Then with probability at least $1-\beta$ where $0<\beta <1$, we have
	\begin{align*}
	\max_{1 \leq k < l \leq K}  \big| \widehat{\mathcal{V}}_{kl} - \mathcal{V}_{kl} \big| \leq 4\sqrt{\frac{B}{n_\text{\emph{min}}}} +  2\sqrt{\frac{B}{n_{\text{\emph{min}}}} \log \bigg\{ \frac{2}{\beta}\binom{K}{2} \bigg\}}. 
	\end{align*}
	\begin{proof}
		Using Theorem 7 of \cite{gretton2012kernel}, one can obtain
		\begin{align*}
		\mP \left( \big| \widehat{\mathcal{V}}_{kl} - \mathcal{V}_{kl} \big| \geq 2 \sqrt{n_k^{-1} B} + 2 \sqrt{n_l^{-1}B} + t \right) \leq 2 \exp \bigg\{ -\frac{(n_k + n_l)\gamma_{k,l}t^2}{2B} \bigg\}.
		\end{align*}
		Then the result follows by applying the union bound as in Theorem~\ref{Theorem: K-sample concentration} and the following inequality 
		\begin{equation*}
		\min_{1 \leq k < l \leq K} (n_k+n_l) \gamma_{k,l} \geq \frac{n_\text{min}}{2}.  \qedhere
		\end{equation*}
	\end{proof}
\end{lemma}

By building on Theorem~\ref{Theorem: K-sample concentration} and Lemma~\ref{Lemma: Unconditional Analysis}, we prove the uniform consistency of the permutation test against $\mathcal{F}_h(\epsilon_N)$ when $\epsilon_N$ is much larger than $r_N^\star$. We provide the proof in Appendix~\ref{Section Appendix}. 

\begin{theorem}[Uniform consistency of the original permutation test] \label{Theorem: Upper Bound}
	Assume that $\textbf{\emph{(B1)}}$ and $\textbf{\emph{(B2)}}$ are fulfilled. Denote the permutation test function by $\phi_{K,\text{\emph{perm}}} = \mathds{1}(p_{\text{\emph{perm}}} \leq \alpha)$ where $p_{\text{\emph{perm}}}$ is given in (\ref{Eq: permutation p-value}). Then under $H_1$, 
	\begin{align*}
	\limsup_{n_{\text{\emph{min}}} \rightarrow \infty} \sup_{(P_1,\ldots,P_K) \in \mathcal{F}_h(b_N r_N^\star)}  \mP\left(\phi_{K,\text{\emph{perm}}}  = 0  \right) = 0,
	\end{align*}
	where $b_N$ is an arbitrary sequence that goes to infinity as $n_{\text{\emph{min}}} \rightarrow \infty$. 
\end{theorem}

Next by using Dvoretzky--Kiefer--Wolfowitz (DKW) inequality \citep[e.g.][]{massart1990tight}, we extend the previous result to the randomized permutation test.

\begin{corollary}[Uniform consistency of the randomized permutation test] \label{Corollary: Upper Bound 2}
	Assume that $\textbf{\emph{(B1)}}$ and $\textbf{\emph{(B2)}}$ are fulfilled. Denote the Monte-Carlo-based permutation test function by $\phi_{K,\text{\emph{MC}}} = \mathds{1}(p_{\text{\emph{MC}}} \leq \alpha)$ where $p_{\text{\emph{MC}}}$ is given in (\ref{Eq: Monte-Carlo p-value}). Then under $H_1$, 
	\begin{align*}
	\lim_{M \rightarrow \infty} \limsup_{n_{\text{\emph{min}}} \rightarrow \infty} \sup_{(P_1,\ldots,P_K) \in \mathcal{F}_h(b_N r_N^\star)}  \mP\left(\phi_{K,\text{\emph{MC}}}  = 0  \right) = 0,
	\end{align*}
	where $b_N$ is an arbitrary sequence that goes to infinity as $n_{\text{\emph{min}}} \rightarrow \infty$. 
\end{corollary}

\begin{remark}
	It is worth pointing out that the results of both Theorem~\ref{Theorem: Upper Bound} and Corollary~\ref{Corollary: Upper Bound 2} hold regardless of whether $K$ is fixed or increases with $n_{\text{min}}$. However, we note that $K$ cannot increase much faster than $e^{n_{\text{min}}}$ as $\max_{1 \leq k < l \leq K} \mathcal{V}_{kl}$ is upper bounded by a positive constant under \textbf{(B1)} and thereby $r_N^\star = \sqrt{\log K / n_{\text{min}}}$ is also bounded.
\end{remark}

\subsection{Minimax rate optimality} \label{Section: Minimax rate optimality}
Theorem~\ref{Theorem: Upper Bound} as well as Corollary~\ref{Corollary: Upper Bound 2} show that the original and randomized permutation tests can be uniformly powerful over $\mathcal{F}_h(b_N r_N^\star)$ when $b_N$ is sufficiently large. In this subsection, we focus on the MMD associated with a translation invariant kernel defined on $\mathbb{R}^d$ and further show that the previous result cannot be improved from a minimax point of view. A kernel $h: \mathbb{R}^d \times \mathbb{R}^d \mapsto \mathbb{R}$ is called \emph{translation invariant} if there exists a symmetric positive definite function $\varphi: \mathbb{R}^d \mapsto \mathbb{R}$ such that $\varphi(x-y) = h(x,y)$ for all $x,y \in \mathbb{R}^d$ \citep{tolstikhin2017minimax}. Then our result is stated as follows. 

\begin{theorem} \label{Theorem: Minimax Lower Bound}
	Let $0<\alpha<1$ and $0 < \zeta < 1-\alpha$. Suppose that $n_{\text{\emph{min}}} \rightarrow \infty$ and $K \rightarrow \infty$. Consider the class of sparse alternatives $\mathcal{F}_h(\epsilon_N)$ defined with a translation invariant kernel $h$ on $\mathbb{R}^d$. Assume that there exists $z \in \mathbb{R}^d$ and $\kappa_1, \kappa_2 > 0$ such that $\varphi(0) - \varphi(z) \geq \kappa_1$ and $r_N^\star \leq \kappa_2$ for all $n_{\text{\emph{min}}}$. Further assume that $\textbf{\emph{(B1)}}$ and $\textbf{\emph{(B2)}}$ hold. Then under $H_1$, there exists a small constant $b > 0$ such that 
	\begin{align*}
	\liminf_{n_{\text{\emph{min}}} \rightarrow \infty} \inf_{\phi \in \Phi_N(\alpha)} \sup_{(P_1,\ldots,P_K) \in \mathcal{F}_h(b r_N^\star)} \mP(\phi = 0) \geq \zeta,
	\end{align*}
	where $\Phi_N(\alpha)$ is the set of all level $\alpha$ test functions such that $\phi: \{Z_1,\ldots,Z_N\} \mapsto \{0,1\}$. 
\end{theorem}

\vskip .5em

\begin{remark}
	The results in Theorem~\ref{Theorem: Upper Bound} and Theorem~\ref{Theorem: Minimax Lower Bound} imply that the proposed permutation test is not only consistent but also minimax rate optimal against the considered sparse alternatives. As far as we are aware, this is the first time that the power of the permutation test is theoretically analyzed under large $N$ and large $K$ situations.
\end{remark}

\begin{remark}
	In our problem setup, a distance between two distributions is measured in terms of the maximum mean discrepancy associated with kernel $h$. One can also study minimax optimality of the proposed test over a class of alternatives measured in terms of a more standard metric such as the $L_2$ distance. For this direction, the results of \cite{li2019optimality} seem useful in which the authors explore minimax rate optimality of kernel mean embedding methods over a Sobolev space in the $L_2$ distance. We leave a detailed analysis of minimax optimality of the proposed test in other metrics to future work.
\end{remark}

\section{Simulations} \label{Section: Simulations}
In this section, we demonstrate the finite-sample performance of the proposed approach via simulations. We consider two characteristic kernels for our test statistic; 1) Gaussian kernel and 2) energy distance kernel. Gaussian kernel is given by $h(x,y) = \exp(-\|x-y\|_2^2/\sigma)$ for which we choose the tuning parameter $\sigma$ by the median heuristic \citep{gretton2012kernel}. On the other hand, energy distance kernel is given by $h(x,y) = (\|x\|_2 + \|y\|_2 - \|x-y\|_2) / 2$ as before. Note that the MMD statistic with energy distance kernel is equivalent to the energy statistic \citep{szekely2004testing,baringhaus2004new} in the two-sample case.

\subsection{Other multivariate $K$-sample tests}
We compare the performance of the proposed tests with two multivariate $K$-sample tests. The first one is the test based on DISCO statistic proposed by \cite{rizzo2010disco}. Let $E_{kl,\alpha^\prime}$ be the $\alpha^\prime$-energy statistic between $P_k$ and $P_l$ given by
\begin{align*}
E_{kl,\alpha^\prime} ~=~ & \frac{2}{n_k n_l} \sum_{i_1=1}^{n_k} \sum_{i_2=1}^{n_l} g_{\alpha^\prime}(X_{i_1,k}, X_{i_2,l})  - \frac{1}{n_k^2}  \sum_{i_1,i_2=1}^{n_k} g_{\alpha^\prime}(X_{i_1,k}, X_{i_2,k}) \\[.5em]
& - \frac{1}{n_l^2}\sum_{i_1,i_2=1}^{n_l} g_{\alpha^\prime}(X_{i_1,l}, X_{i_2,l}),
\end{align*}
where $g_{\alpha^\prime}(x,y) = \|x-y\|_2^{\alpha^\prime}$. Let us write the between-sample and within-sample dispersions by $S_{\alpha^\prime} = K^{-1}\sum_{1 \leq k < l \leq K} E_{kl,\alpha^\prime}$ and $W_{\alpha^\prime} = 2^{-1}\sum_{k=1}^K n_k^{-1}\sum_{i_1, i_2=1}^{n_k} g_{\alpha^\prime}(X_{i_1,k},X_{i_2,k})$. Then DISCO statistic is defined as ratio of the between-sample dispersion to the within-sample dispersion, that is
\begin{align*}
D_\gamma = \frac{S_{\alpha^\prime}/(K-1)}{W_{\alpha^\prime}/(N-K)}.
\end{align*}
The second test, proposed by \cite{huvskova2008tests}, is based on the empirical characteristic functions. For a given $\alpha^{\prime\prime} \in \mathbb{R}$,  \cite{huvskova2008tests} consider the weighted $L_2$ distance between empirical characteristic functions as their test statistic, that is
\begin{align*}
H_{\alpha^{\prime\prime}} ~ = ~  \sum_{k=1}^K \frac{N-n_k}{Nn_k} \sum_{i_1, i_2=1}^{n_k} e^{-\|X_{i_1,k} - X_{i_2,k} \|_2^2/4\alpha^{\prime\prime}} - \frac{1}{N} \sum_{1 \leq k \neq l \leq K} \sum_{i_1=1}^{n_k} \sum_{i_2=1}^{n_l} e^{-\|X_{i_1,k} - X_{i_2,l} \|_2^2/4\alpha^{\prime\prime}}.
\end{align*}
In their paper, \cite{huvskova2008tests} consider $\alpha^{\prime\prime} = 1,1.5,2$ in their simulation study. Throughout our simulations, we choose $\alpha^\prime=1$ for $D_{\alpha^\prime}$ and $\alpha^{\prime\prime} =1.5$ for $H_{\alpha^{\prime\prime}}$ and reject the null for large values of $D_{\alpha^\prime}$ and $H_{\alpha^{\prime\prime}}$. 

We also attempted to consider the graph-based $K$-sample test recently developed by \cite{wang2018nonparametric}. To implement their test, we used the R package provided by the same authors. Unfortunately, their method was not applicable when $K$ is large due to numerical overflow in computing orthogonal polynomials. Hence we focus on the first two methods described in this subsection and compare them with the proposed tests against sparse alternatives.

\subsection{Set-up}
Let us denote a multivariate normal distribution with mean vector $\mu$ and covariance matrix $\Sigma$ by $N(\mu, \Sigma)$. Similarly we denote a multivariate Laplace distribution with mean vector $\mu$ and covariance matrix $\Sigma$ by $L(\mu,\Sigma)$. We examine the performance of the considered tests under the following sparse alternatives:
\begin{itemize}
	\item[(a)] {\textbf{Normal Location}}: $P_1 = N(\delta_1, I_d)$ and $P_2=\ldots=P_K=N(\delta_0, I_d)$, \\[-.5em]
	\item[(b)] {\textbf{Normal Scale}}: $P_1 = N(\delta_0, 3 \times I_d)$ and $P_2=\ldots=P_K=N(\delta_0, I_d)$, \\[-.5em]
	\item[(c)] {\textbf{Laplace Location}}: $P_1 = L(\delta_{1.2}, I_d)$ and $P_2=\ldots=P_K=L(\delta_0, I_d)$, \\[-.5em]
	\item[(d)] {\textbf{Laplace Scale}}: $P_1 = L(\delta_0, 3 \times I_d)$ and $P_2 = \ldots=P_K=L(\delta_0, I_d)$,
\end{itemize}
where $\delta_b = (b,\ldots,b)^\top$ and $I_d$ is the $d$-dimensional identity matrix. In words, we consider the sparse alternatives where only one of the distributions differs from the other $K-1$ distributions. Consequently, the signal is getting sparser as $K$ increases. Throughout our experiments, we fix sample sizes $n_1=n_2=\ldots=n_K=10$ and dimension $d=5$ while increasing the number of distributions $K\in \{2,20,40,60,80,100\}$. All tests were implemented via the randomized permutation procedure with $M=200$ random permutations using the $p$-value in (\ref{Eq: Monte-Carlo p-value}). As a result, they are all valid level $\alpha$ tests. Simulations were repeated $800$ times to estimate the power at significance level $\alpha = 0.05$. 

\begin{figure}[!t]
	\centering
	\subfigure[]{\includegraphics[width=0.48\textwidth]{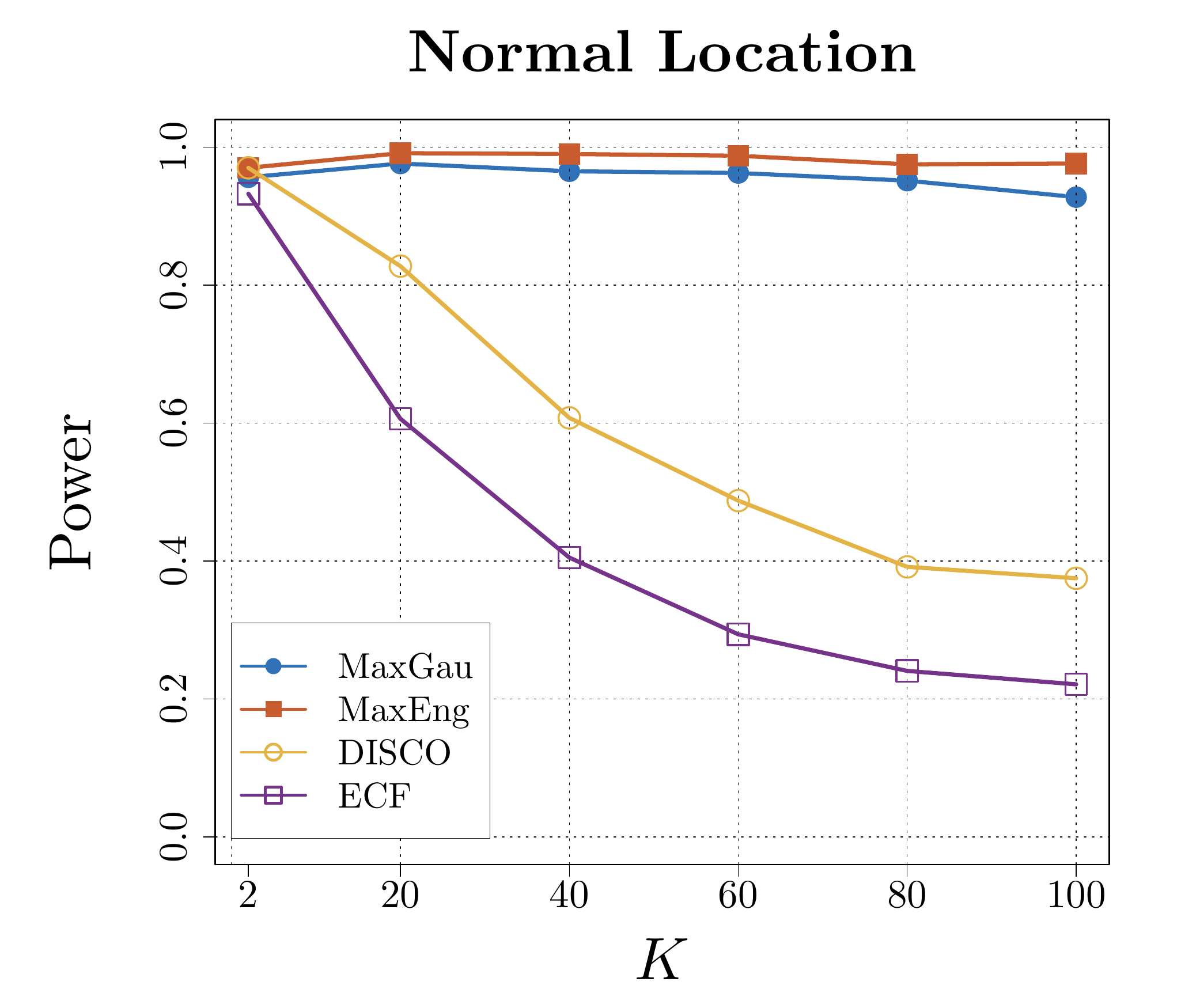}} 
	\subfigure[]{\includegraphics[width=0.48\textwidth]{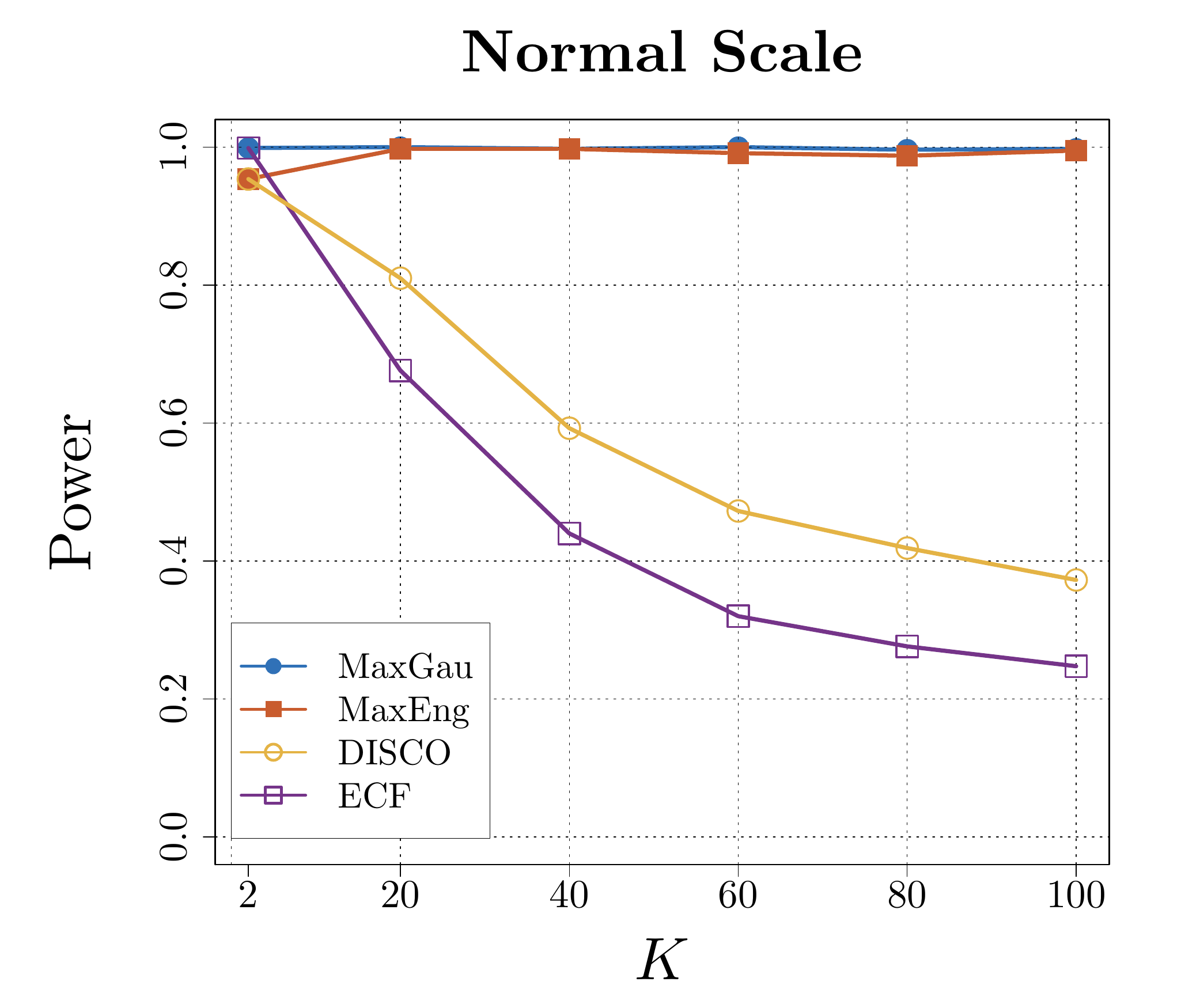}} 
	\subfigure[]{\includegraphics[width=0.48\textwidth]{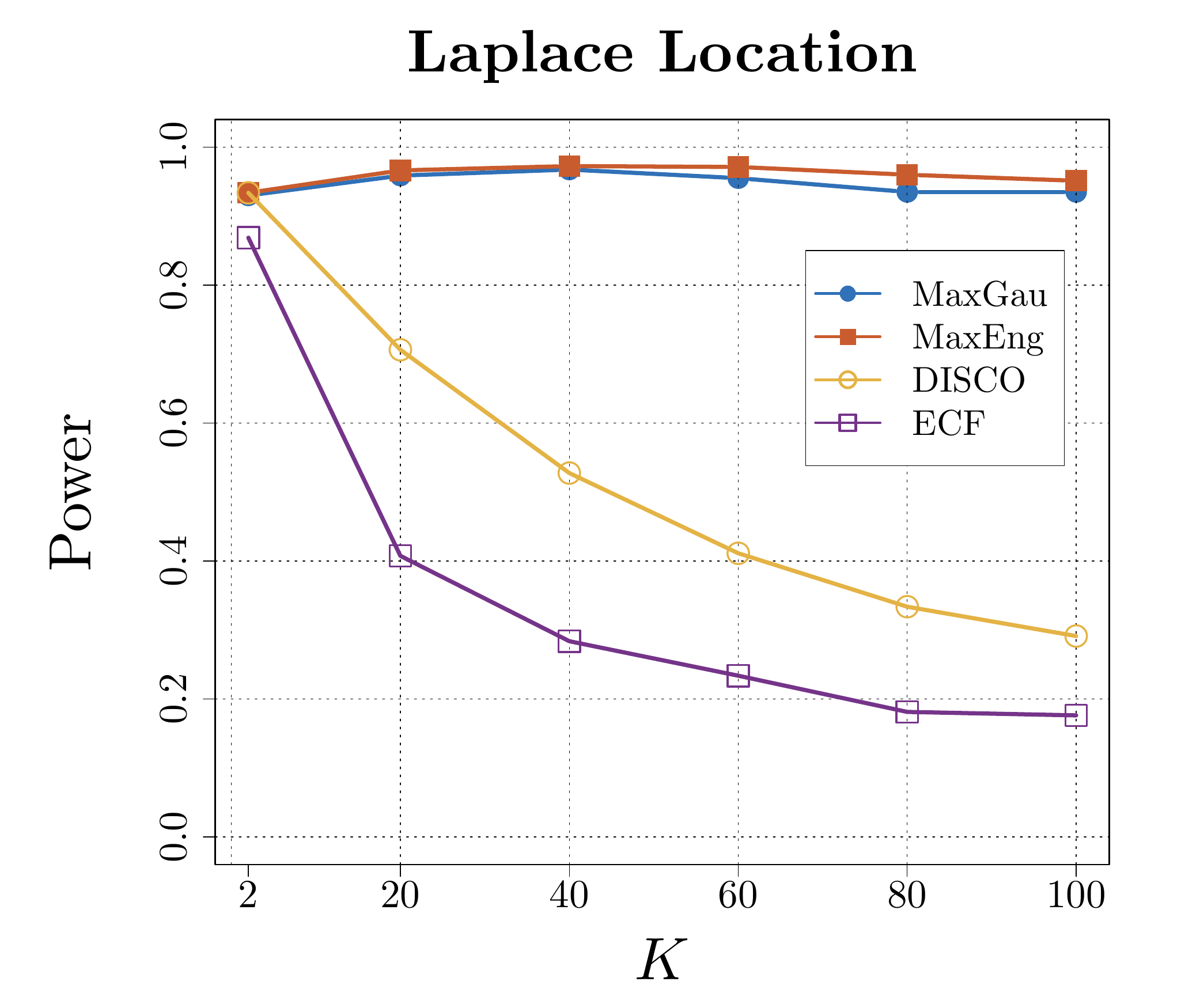}}
	\subfigure[]{\includegraphics[width=0.48\textwidth]{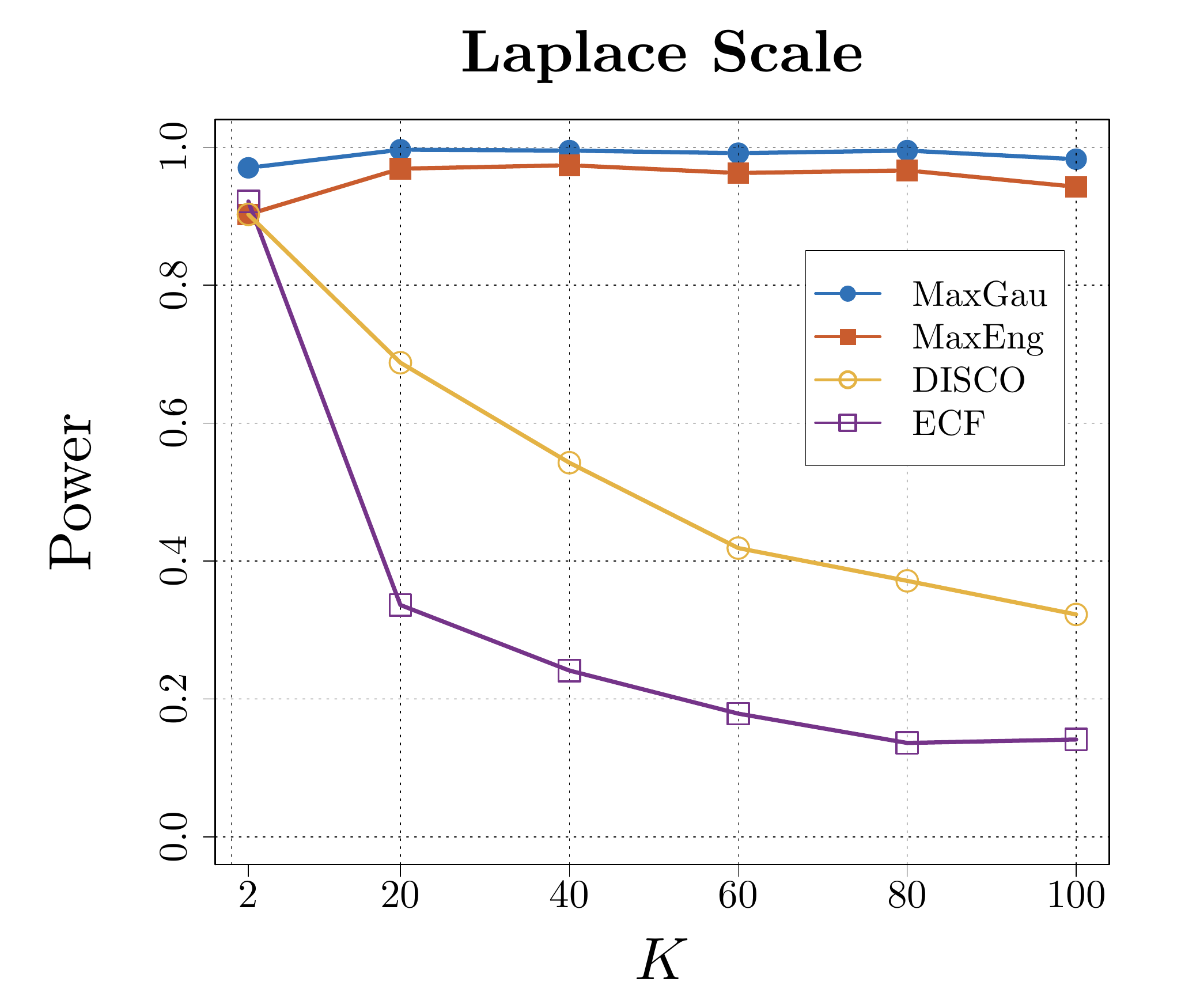}}
	\caption{Empirical power comparisons of the considered tests against (a)~Normal location, (b)~Normal scale, (c)~Laplace location, (d)~Laplace scale alternatives. We refer to the tests based on $\widehat{\mathcal{V}}_{h,\text{max}}$ with Gaussian kernel and energy distance kernel as MaxGau and MaxEng, respectively. In addition, the tests based on $D_{\alpha^\prime}$ and $H_{\alpha^{\prime\prime}}$ are referred to as DISCO and ECF, respectively. See Section~\ref{Section: Simulations} for details.}
	\label{Figure: Power}
\end{figure}

\subsection{Results}
From the results presented in Figure~\ref{Figure: Power}, we observe that the tests based on $D_{\alpha^\prime}$ and $H_{\alpha^{\prime\prime}}$ have consistently decreasing power as $K$ increases in all sparse scenarios. This  can be explained by the fact that $D_{\alpha^\prime}$ and $H_{\alpha^{\prime\prime}}$ are defined as an average between pairwise distances. Under the given sparse scenario, the average of pairwise distances, which is a signal to reject $H_0$, decreases as $K$ increases. Hence the resulting tests based on $D_{\alpha^\prime}$ and $H_{\alpha^{\prime\prime}}$ suffer from low power in large $K$. On the other hand, the proposed tests show robust performance to the number of distributions $K$ under the given setting. They in fact have power very close to one even when $K$ is considerably large, which emphasizes the benefit of using the maximum-type statistic against sparse alternatives. 

Despite their good performance over sparse alternatives, the proposed tests do not always perform better than the average-type tests based on $D_{\alpha^\prime}$ and $H_{\alpha^{\prime\prime}}$. For example, these average-type tests may outperform the proposed maximum-type tests against dense alternatives where many of $P_1,\ldots,P_K$ differ from each other. Given that prior knowledge on alternatives is not always available to users, developing a powerful test against both dense and sparse alternatives is an interesting direction for future work.


\section{Conclusion} \label{Section: Conclusion}
In this paper, we introduced a new nonparametric $K$-sample test based on the maximum mean discrepancy. The limiting distribution of the proposed test statistic was derived based on Cram{\'e}r-type moderate deviation for degenerate two-sample $V$-statistics. Unfortunately, the limiting distribution relies on an infinite number of nuisance parameters, which are intractable in general. Due to this challenge, we considered the permutation approach to determine the cut-off value of the test. 
We provided a concentration inequality for the proposed test statistic with a sharp exponential tail bound under permutations. On the basis of this result, we studied the power of the permutation test in large $K$ and large $N$ situations and further proved its minimax rate optimality under some regularity conditions. From our simulation studies, the proposed test is shown to be powerful against sparse alternatives where the previous methods suffer from low power. These findings suggest that our method will be useful in application areas where only a small number of populations differ from the others. 

The power analysis in Section~\ref{Section: Power Analysis} relies on the assumption that a kernel is uniformly bounded. Although some of the popular kernels satisfy this assumption, our result cannot be applied to unbounded cases. One possible way to address this issue is to impose appropriate moment conditions on a kernel and utilize a suitable concentration inequality  \citep[e.g. a modified McDiarmid's inequality in][]{kontorovich2014concentration} to obtain a similar result to Lemma~\ref{Lemma: Unconditional Analysis}. This topic is reserved for future work. 

\section*{Acknowledgements}
The author would like to thank the associate editor and the anonymous reviewers for their valuable comments. The author also thanks Sivaraman Balakrishnan and Larry Wasserman for their kind support and constructive feedback.

\bibliographystyle{apalike}
\bibliography{reference}

\clearpage

\appendix

\allowdisplaybreaks

\section{Appendix} \label{Section Appendix}

In this section, we collect the proofs of the theorems in the main text. Throughout this section, we use $C_1,C_2,\ldots$ to denote some constants that may change from line to line.

\subsection{Proof of Theorem~\ref{Theorem: moderate deviation}}
The following proof is built upon the proof of Theorem 4.1 of \cite{drton2018high} and extends theirs to two-sample $V$-statistics and unbounded eigenfunctions. We start with another representation of $\widehat{\mathcal{V}}_{12}^2$ in terms of $\{\lambda_v\}_{v=1}^\infty$ and $\{ \varphi_v(\cdot) \}_{v=1}^\infty$. Since $h(z_1,z_2)$ is symmetric in its arguments, $\widehat{\mathcal{V}}_{12}^2$ can also be represented in terms of the centered kernel as
\begin{align*}
\widehat{\mathcal{V}}_{12}^2 = \frac{1}{n_1^2}\sum_{i_1,i_2=1}^{n_1} \overline{h}(X_{i_1,1}, X_{i_2,1}) + \frac{1}{n_2^2}\sum_{i_1,i_2=1}^{n_2} \overline{h}(X_{i_1,2}, X_{i_2,2})- \frac{2}{n_1 n_2} \sum_{i_1=1}^{n_1} \sum_{i_2=1}^{n_2} \overline{h}(X_{i_1,1}, X_{i_2,2}).
\end{align*}
Furthermore, based on the decomposition given in (\ref{Eq: centered kernel}), $\widehat{\mathcal{V}}_{12}^2$ can be written as
\begin{align*}
\widehat{\mathcal{V}}_{12}^2 ~=~ \sum_{v=1}^\infty \lambda_v \Bigg\{ \frac{1}{n_1} \sum_{i_1=1}^{n_1} \varphi_v (X_{i_1,1}) -  \frac{1}{n_2} \sum_{i_2=1}^{n_2} \varphi_v (X_{i_2,2}) \Bigg\}^2.
\end{align*}
In what follows, we consider two different cases: 1) $x$ is bounded and 2) $x$ tends to infinity and prove Theorem~\ref{Theorem: moderate deviation} under each scenario. 

\vskip 1em

\noindent \textbf{Case 1: $x$ is bounded}

\vskip .5em
\noindent First write the corresponding degenerate two-sample $U$-statistic by 
\begin{align*}
\widehat{\mathcal{U}}_{12} ~ = ~ & \frac{1}{n_1(n_1-1)}\sum_{1 \leq i_1 \neq i_2 \leq n_1} \overline{h}(X_{i_1,1}, X_{i_2,1}) + \frac{1}{n_2(n_2-1)}\sum_{1\leq i_1 \neq i_2 \leq n_2} \overline{h}(X_{i_1,2}, X_{i_2,2}) \\[.5em]
& - \frac{2}{n_1 n_2} \sum_{i_1=1}^{n_1} \sum_{i_2=1}^{n_2} \overline{h}(X_{i_1,1}, X_{i_2,2}).
\end{align*}
Then using the result on Chapter 3 of \cite{bhat1995theory}, 
\begin{align*}
\frac{n_1n_2}{N}\widehat{\mathcal{U}}_{12} ~ \convD ~ \sum_{v=1}^{\infty} \lambda_v (\xi_v^2 - 1). 
\end{align*}
Now the difference between the $V$-statistic and $U$-statistic is 
\begin{align*}
\widehat{\mathcal{V}}_{12}^2 - \widehat{\mathcal{U}}_{12} ~=~ &\frac{1}{n_1^2} \sum_{i_1=1}^{n_1} \overline{h}(X_{i_1,1}, X_{i_1,1}) +  \frac{1}{n_2^2} \sum_{i_2=1}^{n_2} \overline{h}(X_{i_2,2}, X_{i_2,2}) \\[.5em]
& - \frac{1}{n_1^2(n_1-1)} \sum_{1 \leq i_1 \neq i_2 \leq n_1} \overline{h}(X_{i_1,1}, X_{i_2,1}) - \frac{1}{n_2^2(n_2-1)} \sum_{1 \leq i_1 \neq i_2 \leq n_2} \overline{h}(X_{i_1,2}, X_{i_2,2}). 
\end{align*}
Under the assumption that $\mE[|\overline{h}(X_1,X_1)|] < \infty$, we apply the strong law of large numbers for $U$-statistics \citep[e.g. Theorem A of Section 5.4 in][]{serfling1980approximation} to have
\begin{align*}
\frac{n_1n_2}{N}\left( \widehat{\mathcal{V}}_{12}^2 - \widehat{\mathcal{U}}_{12} \right) ~ \convAS ~ \mE[\overline{h}(X_1,X_1)] = \sum_{v=1}^\infty \lambda_v.
\end{align*}
Hence we establish that 
\begin{align*}
\frac{n_1n_2}{N}\widehat{\mathcal{V}}_{12}^2 ~ \convD ~ \sum_{v=1}^{\infty} \lambda_v \xi_v^2,
\end{align*}
which leads to (\ref{Eq: Cramer-type moderate deviation}) for any bounded $x$.

\vskip 1.5em

\noindent \textbf{Case 2: $x$ tends to infinity}

\vskip .5em

\noindent Next we focus on the case where $x$ tends to infinity at a certain rate. To start, for a sufficiently large positive integer $T$ to be specified later, let us define the truncated statistic
\begin{align*}
\widehat{\mathcal{V}}^2_T ~=~ \sum_{v=1}^T \lambda_v \Bigg\{ \frac{1}{\sqrt{n_1}} \sum_{i=1}^{n_1} \varphi_v (X_{i,1}) -  \frac{1}{\sqrt{n_2}} \sum_{i=1}^{n_2} \varphi_v (X_{i,2}) \Bigg\}^2.
\end{align*}
Based on Slutsky's argument,
\begin{align*}
\mP\left( \frac{n_1n_2}{N} \widehat{\mathcal{V}}_{12}^2 \geq x \right) & ~ \leq ~ \mP\left( \frac{n_1n_2}{N} \widehat{\mathcal{V}}_T^2 \geq x  - \epsilon_1 \right) + \mP \bigg\{ \Big|\frac{n_1n_2}{N} \left(\widehat{\mathcal{V}}_{12}^2  -  \widehat{\mathcal{V}}_T^2 \right) \Big| \geq \epsilon_1 \bigg\} \\[.5em]
& ~ := ~ (I) + (II) \quad \text{(say)}.
\end{align*}
Here and hereafter $\epsilon_1, \epsilon_2,\epsilon_3$ are some positive constants that will be specified later. Let us rewrite
\begin{align*}
\sqrt{\frac{\lambda_v}{n_1}} \sum_{i_1=1}^{n_1} \varphi_v (X_{i_1,1}) -   \sqrt{\frac{\lambda_v}{n_2}}  \sum_{i_2=1}^{n_2} \varphi_v (X_{i_2,2})  = \sum_{i=1}^{N} \sqrt{\lambda_v} w_i \varphi_v(Z_i)
\end{align*}
where
\begin{align*}
& (w_1,\ldots,w_N) ~=~ (n_1^{-1/2},\ldots,n_1^{-1/2},-n_2^{-1/2},\ldots,-n_2^{-1/2}),\\[.5em]
& (Z_1,\ldots,Z_N) ~=~ (X_{1,1},\ldots,X_{n_1,1},X_{1,2},\ldots,X_{n_2,2}).
\end{align*}
Further let $\varphi_{1,\ldots,T}^{\lambda}(Z_i) = (\sqrt{\lambda_1}w_i\varphi_1(Z_i),\ldots,\sqrt{\lambda_T}w_i \varphi_T(Z_i))^\top$. For each $i=1,\ldots,N$, we verify the multivariate Bernstein condition used in \cite{zaitsev1987gaussian}. Specifically, for any $u,v \in \mathbb{R}^{T}$ and $m=3,4,\ldots$, we have that
\begin{align*}
&\big| \mE \big[ \{ \varphi_{1 \cdots T}^{\lambda}(Z_i)^\top u \}^2 \{\varphi_{1 \cdots T}^{\lambda}(Z_i)^\top v\}^{m-2} \big] \big| \\[.5em]
~ \overset{(i)}{\leq} ~ &  \left( \sqrt{\frac{\lambda_1}{n_1}} +  \sqrt{\frac{\lambda_1}{n_2}}  \right)^m \big| \mE \big[ \{ \varphi_{1 \cdots T}(Z_i)^\top u \}^2 \{\varphi_{1 \cdots T}(Z_i)^\top v \}^{m-2} \big] \big| \\[.5em]
~ \overset{(ii)}{\leq} ~ & \left( \sqrt{\frac{\lambda_1}{n_1}} +  \sqrt{\frac{\lambda_1}{n_2}}  \right)^m \gamma^m m^{m/2} \|u\|_2^2 \|v\|_2^{m-2} \\[.5em]
~ \overset{(iii)}{\leq} ~ & \left( \sqrt{\frac{\lambda_1}{n_1}} +  \sqrt{\frac{\lambda_1}{n_2}}  \right)^m \gamma^m m!\|v\|_2^{m-2}  \mE\big[ \{ \varphi_{1 \cdots T}(Z_i)^\top u \}^2\big] 
\end{align*}
where 
\begin{itemize}
	\item $(i)$ follows since $\lambda_1 \leq \lambda_2 \leq \ldots$ and $\max\{1/\sqrt{n_1}, 1/\sqrt{n_2}\} \leq 1/\sqrt{n_1} + 1/\sqrt{n_2}$. 
	\item $(ii)$ uses the condition \textbf{(A2)}.
	\item $(iii)$ uses $m! \geq m^{m/2}$ for all $m \geq 3$ and $\mE\big[ \{ \varphi_{1 \cdots T}(Z_i)^\top u \}^2\big]  = \| u \|_2^2$. 
\end{itemize}
Thus together with the assumption that $C_1^{-1} \leq n_1/n_2 \leq C_1$, the multivariate Bernstein condition in \cite{zaitsev1987gaussian} is fulfilled with his notation $\tau = C_2 N^{-1/2}$ for sufficiently large $C_2$. Consequently, we can apply Theorem 1.1 of \cite{zaitsev1987gaussian} to show that
\begin{align*}
\mP\left( \frac{n_1n_2}{N} \widehat{\mathcal{V}}_T^2 \geq x  - \epsilon_1 \right) ~ \leq ~ & \mP \left[ \sum_{v=1}^T \lambda_v \xi_v^2 \geq \{(x - \epsilon_1)^{1/2} - \epsilon_2 \}^2 \right] \\[.5em]
& + C_{3} T^{5/2} \exp\left( - \frac{\sqrt{N}\epsilon_2}{C_{4} T^{5/2}} \right).
\end{align*}
By applying Slutsky's argument again, the first term is bounded by
\begin{align*}
\mP \left[  \sum_{v=1}^T \lambda_v  \xi_v^2 \geq \{(x - \epsilon_1)^{1/2} - \epsilon_2 \}^2 \right] ~ \leq ~ & \mP \left[ \sum_{v=1}^\infty \lambda_v  \xi_v^2 \geq \{(x - \epsilon_1)^{1/2} - \epsilon_2 \}^2 - \epsilon_3 \right] \\[.5em]
& + \mP \left(  \bigg| \sum_{v=T+1}^\infty \lambda_v\xi_v^2 \bigg| \geq  \epsilon_3 \right). 
\end{align*}

For a random variable $X$, let us denote the sub-Gaussian norm and sub-exponential norm by $\| X \|_{\psi_2} := \inf \{t > 0: \mE[\exp(X^2/t^2)] \leq 2 \}$ and $\| X \|_{\psi_1} := \inf \{t > 0: \mE[\exp(|X|/t)] \leq 2 \}$, respectively. By the property of the norm, Example 2.5.8 of \cite{vershynin2018high} and Lemma 2.7.6 of \cite{vershynin2018high}, we observe that 
\begin{align*}
\bigg\| \sum_{v=T+1}^\infty \lambda_v \xi_v^2 \bigg\|_{\psi_1} \leq \sum_{v=T+1}^\infty \lambda_v \big\|  \xi_v^2 \big\|_{\psi_1} =  \sum_{v=T+1}^\infty \lambda_v \big\|  \xi_v \big\|_{\psi_2}^2 \leq C_{5} \sum_{v=T+1}^\infty \lambda_v.
\end{align*}
Then by Proposition 2.7.1 of \cite{vershynin2018high},
\begin{align*}
\mP \left( \bigg| \sum_{v=T+1}^\infty \lambda_v \xi_v^2 \bigg|   \geq \epsilon_3  \right) \leq 2 \exp \left( - \frac{\epsilon_3}{C_{5} \sum_{v=T+1}^\infty \lambda_v} \right).
\end{align*}
Thus
\begin{align*}
(I) ~ \leq  ~ & \mP \left[  \sum_{v=1}^\infty \lambda_v  \xi_v^2 \geq \{(x - \epsilon_1)^{1/2} - \epsilon_2 \}^2 - \epsilon_3 \right]   \\[.5em]
& +   C_{3} T^{5/2} \exp\left( - \frac{\sqrt{N}\epsilon_2}{C_{4} T^{5/2}} \right) + 2 \exp \left( - \frac{\epsilon_3}{C_{5} \sum_{v=T+1}^\infty \lambda_v} \right).
\end{align*}

Next we focus on the term $(II)$. Note that the multivariate moment condition in (\ref{Eq: Multivariate moment condition}) implies the univariate sub-Gaussian condition for $\varphi_v(X_1)$ and $v=1,2,\ldots$. That is, there exists a constant $C_{6} >0$ independent of $v$ such that 
\begin{align*}
\mE[\{\varphi_v(X_1)\}^m] \leq C_{6} m^{m/2} \mE[\{\varphi_v(X_1)\}^2] = C_{6} m^{m/2} \quad \text{for all $m \geq 1$.}
\end{align*}
Thus, followed by Proposition 2.7.1 of \cite{vershynin2018high}, $\varphi_v(X_1)$ has a finite sub-Gaussian norm and furthermore $\sup_{v \geq 1} \| \varphi_v(X_1) \|_{\psi_2} := C_{7} < \infty$. Then
\begin{align*}
\bigg\|  \frac{n_1n_2}{N} \left(\widehat{\mathcal{V}}^2  -  \widehat{\mathcal{V}}_T^2 \right)  \bigg\|_{\psi_{1}} ~ \overset{(i)}{\leq} ~ &  \frac{n_1n_2}{N}\sum_{v=T+1}^\infty \lambda_v \Bigg\| \Bigg[ \frac{1}{n_1} \sum_{i_1=1}^{n_1} \varphi_v(X_{i_1,1})  - \frac{1}{n_2} \sum_{i_2=1}^{n_2} \varphi_v(X_{i_2,2})\Bigg]^2 \Bigg\|_{\psi_{1}} \\[.5em]
~ \overset{(ii)}{=}  ~ &  \frac{n_1n_2}{N}\sum_{v=T+1}^\infty \lambda_v \Bigg\|\frac{1}{n_1} \sum_{i_1=1}^{n_1} \varphi_v(X_{i_1,1})  - \frac{1}{n_2} \sum_{i_2=1}^{n_2} \varphi_v(X_{i_2,2}) \Bigg\|_{\psi_{2}}^2 \\[.5em]
~ \overset{(iii)}{\leq}  ~& C_{8} \frac{n_1n_2}{N}\sum_{v=T+1}^\infty \lambda_v \Bigg[\frac{1}{n_1} \sum_{i_1=1}^{n_1} \big\| \varphi_v(X_{i_1,1}) \big\|_{\psi_{2}}^2  + \frac{1}{n_2} \sum_{i_2=1}^{n_2} \big\| \varphi_v(X_{i_2,2}) \big\|_{\psi_{2}}^2  \Bigg] \\[.5em]
~ \overset{(iv)}{\leq}  ~&  C_{9} \sum_{v=T+1}^\infty \lambda_v
\end{align*}
where 
\begin{itemize}
	\item $(i)$ uses the triangle inequality.
	\item $(ii)$ uses Lemma 2.7.6 of \cite{vershynin2018high}.
	\item $(iii)$ holds by Proposition 2.6.1 of \cite{vershynin2018high}.
	\item $(iv)$ follows since $\sup_{v \geq 1} \| \varphi_v(X_1) \|_{\psi_2} < \infty$.
\end{itemize}
Based on the above result, we apply Markov's inequality to bound
\begin{align*}
(II) ~ \leq ~ \exp\left( -\frac{\epsilon_1}{C_{9} \sum_{v=T+1}^\infty \lambda_v} \right).
\end{align*}
To summarize, we have obtain that 
\begin{align} \nonumber
\frac{\mP( n_1n_2\widehat{\mathcal{V}}^2 / N \geq x  )}{\mP (  \sum_{v=1}^\infty \lambda_v  \xi_v^2 \geq x )} ~ \leq ~ &\Bigg\{ \mP \left(  \sum_{v=1}^\infty \lambda_v  \xi_v^2 \geq x\right)\Bigg\}^{-1} \cdot \Bigg\{ \mP \left(  \sum_{v=1}^\infty \lambda_v  \xi_v^2 \geq \{(x - \epsilon_1)^{1/2} - \epsilon_2 \}^2 - \epsilon_3 \right)   \\[.5em] \nonumber
& +   C_{3} T^{5/2} \exp\left( - \frac{\sqrt{N}\epsilon_2}{C_{4} T^{5/2}} \right) + 2 \exp \left( - \frac{\epsilon_3}{C_{5} \sum_{v=K+1}^\infty \lambda_v} \right) \\[.5em]
& +  \exp\left( -\frac{\epsilon_1}{C_{9} \sum_{v=T+1}^\infty \lambda_v} \right) \Bigg\}.  \label{Eq: summary}
\end{align}
Our goal is now to show that the right-hand side of (\ref{Eq: summary}) converges to one by properly choosing $x,\epsilon_1,\epsilon_2,\epsilon_3,T$. To simplify the notation, we let $\zeta \overset{d}{=} \sum_{v=1}^\infty \lambda_v \xi_v^2$ and denote $\overline{F}_\zeta(x)= \mP (\sum_{v=1}^\infty \lambda_v  \xi_v^2 \geq x)$. We also write the density function of $\zeta$ by $f_\zeta(x)$.

We start with the first term of (\ref{Eq: summary}). Write
\begin{align*}
\epsilon := x -  \big\{ (x - \epsilon_1)^{1/2} - \epsilon_2 \big\}^2 -\epsilon_3.
\end{align*}
Followed by \cite{zolotarev1961concerning}, we can approximate the survival function and density function of $\zeta$ as
\begin{align*}
& \overline{F}_{\zeta}(x)= \frac{\kappa}{\Gamma(\mu_1/2)} \left( \frac{x}{2\lambda_1} \right)^{\mu_1/2-1} \exp\left( - \frac{x}{2\lambda_1} \right) \{1 + o(1)\} \\[.5em]
& f_{\zeta}(x) = \frac{\kappa}{2\lambda_1 \Gamma(\mu_1/2)} \left( \frac{x}{2\lambda_1} \right)^{\mu_1/2-1} \exp\left( - \frac{x}{2\lambda_1} \right) \{1 + o(1)\}
\end{align*}
for all $x > - \sum_{v=1}^\infty \lambda_v := - \Lambda$ that tends to infinity and $\kappa = \prod_{v=\mu_1+ 1}^\infty(1 - \lambda_v/\lambda_1)^{-1/2}$. Then followed similarly by (A.13) of \cite{drton2018high}, it is seen that there exists a constant $x_0$ such that for all $0 < \epsilon \leq \lambda_1/2$,
\begin{align*}
\sup_{x \geq x_0} \big| \overline{F}^{-1}_{\zeta}(x) \cdot \max_{x^\prime \in [x-\epsilon,x]} f_{\zeta}(x^\prime) \big| \leq 2\lambda_1^{-1} .
\end{align*}  
Using this, the first term is bounded by
\begin{align*}
\frac{\mP \left[  \sum_{v=1}^\infty \lambda_v  \xi_v^2 \geq \{(x - \epsilon_1)^{1/2} - \epsilon_2 \}^2 - \epsilon_3 \right]}{\overline{F}_{\zeta}(x)} ~ \leq ~ &\frac{\mP \left(  \sum_{v=1}^\infty \lambda_v  \xi_v^2 \geq x\right)}{\overline{F}_{\zeta}(x)} + \frac{\epsilon \cdot \max_{x^\prime \in [x-\epsilon,x]} f_{\zeta}(x^\prime)}{\overline{F}_{\zeta}(x)} \\[.5em]
~\leq ~ & 1 + 2\epsilon \lambda_1^{-1}
\end{align*}
for all $x \geq x_0$. Next we shall choose $\epsilon_1,\epsilon_2,\epsilon_3$ decreasing in $N$ so that $1 + 2\epsilon\lambda_1^{-1}$ converges uniformly to one for all $x \geq x_0$. Thus the upper bound of the first term converges to one uniformly over $x \geq x_0$. Hence, we only need to study the last three terms in (\ref{Eq: summary}) to finish the proof.

Let us first specify $T = \floor{N^{(1-3\theta)/5}}$ where $\theta$ satisfies
\begin{align*}
\theta < \sup \Big\{ q \in [0,1/3): \sum\nolimits_{v > \floor{N^{(1-3q)/5}}} \lambda_v = O(N^{-q}) \Big\}.
\end{align*}
Note that by the definition of $\theta$, there exists a positive constant $C_{10}$ such that $\sum_{v=T+1}^\infty \lambda_v \leq C_{10} N^{-\theta}$  for a sufficiently large $N$. Hence it now suffices to show that for all $x \in (0, o(N^{\theta}))$,
\begin{equation}
\begin{aligned} \label{Eq: three inequalities}
& \Bigg\{ \left( \frac{x}{2\lambda_1} \right)^{\mu_1/2-1} \exp\left( - \frac{x}{2\lambda_1} \right)  \Bigg\}^{-1} \exp \left( -\frac{\epsilon_1}{C_{11} N^{-\theta}} \right)~ \leq ~ o(1), \\[.5em]
& \Bigg\{ \left( \frac{x}{2\lambda_1} \right)^{\mu_1/2-1} \exp\left( - \frac{x}{2\lambda_1} \right)  \Bigg\}^{-1} N^{(1-3\theta)/2} \exp\left( - \frac{\sqrt{N}\epsilon_2}{C_{4} N^{(1-3\theta)/2}} \right) 
~ \leq ~ o(1), \\[.5em]
& \Bigg\{ \left( \frac{x}{2\lambda_1} \right)^{\mu_1/2-1} \exp\left( - \frac{x}{2\lambda_1} \right)  \Bigg\}^{-1}   \exp \left( - \frac{\epsilon_3}{C_{12} N^{-\theta}} \right) ~ \leq ~ o(1).
\end{aligned}
\end{equation}
For this purpose, we choose $\epsilon_1$, $\epsilon_2$ and $\epsilon_3$ such that
\begin{align*}
\epsilon_1 = C_{11} N^{-\theta} \left( \frac{x}{2\lambda_1}  + N^{\theta/2} \right), \quad \epsilon_2 = N^{-\theta/2}, \quad \epsilon_3 =  C_{12} N^{-\theta} \left( \frac{x}{2\lambda_1}  + N^{\theta/2} \right),
\end{align*}
which tend to zero as $N \rightarrow \infty$ under $x \in (0, o(N^{\theta}))$. It is then straightforward to see that the three inequalities in (\ref{Eq: three inequalities}) hold under the given setting. Consequently, 
\begin{align*}
\frac{\mP ( n_1n_2\widehat{\mathcal{V}}_{12}^2 / N \geq x  )}{\mP \left(  \sum_{v=1}^\infty \lambda_v  \xi_v^2 \geq x\right)} ~ \leq ~ 1 + o(1).
\end{align*}
The other direction follows similarly, which concludes
\begin{align*}
\frac{\mP ( n_1n_2\widehat{\mathcal{V}}_{12}^2 / N \geq x )}{\mP \left(  \sum_{v=1}^\infty \lambda_v  \xi_v^2 \geq x\right)} ~ = ~ 1 + o(1)
\end{align*}
uniformly over $x \in (0, o(N^{\theta}))$. This completes the proof. 

\vskip 1em

\subsection{Proof of Theorem~\ref{Theorem: Gumbel Limit}} 
Continuing our discussion from Section~\ref{Section: Gumbel limiting distribution}, we apply Lemma~\ref{Lemma: Arratia} together with Theorem~\ref{Theorem: moderate deviation} to obtain the result. Specifically, we set 
\begin{align*}
x = 4\lambda_1 \log K  + \lambda_1 (\mu_1 - 2) \log \log K + \lambda_1 y.
\end{align*}
Then by the triangle inequality
\begin{align*}
& \bigg|\mP \big(n\widehat{ \mathcal{V}}_{h,\text{max}}^2/2 \leq x \big) - \exp\bigg\{ - \frac{2^{\mu_1/2-2}\kappa}{\Gamma(\mu_1/2)}  \exp\left( -\frac{y}{2}\right) \bigg\}  \bigg| \\[.5em]
\leq  ~ & \bigg|\mP \big(n\widehat{ \mathcal{V}}_{h,\text{max}}^2/2 \leq x \big) - \exp\bigg\{ - \frac{K(K-1)}{2} \mP \big( n \widehat{\mathcal{V}}_{12}^2  / 2 > x \big) \bigg\}\bigg| \\[.5em]
& + \bigg| \exp\bigg\{ - \frac{2^{\mu_1/2-2}\kappa}{\Gamma(\mu_1/2)}  \exp\left( -\frac{y}{2}\right) \bigg\} - \exp\bigg\{ - \frac{K(K-1)}{2} \mP \big( n \widehat{\mathcal{V}}_{12}^2  / 2 > x \big) \bigg\}\bigg| \\[.5em]
= ~ & (I) + (II) \quad \text{(say).}
\end{align*}
By setting $\mathcal{I} = \{(i,j): 1 \leq i < j \leq K  \}$ and $B_{u_{i,j}} = \{(k,l) \in \mathcal{I}: \text{Card}\{(k,l) \cap (i,j) \} \neq 0 \}$ where $u_{i,j} := (i,j) \in \mathcal{I}$ and $\text{Card}\{A\}$ denotes the cardinality of a set $A$, Lemma~\ref{Lemma: Arratia} yields  $(I) \leq b_1 + b_2 + b_3$. Here, in our setting,
\begin{align*}
& b_1 = \frac{K(K-1)(2K-3)}{2} \big\{ \mP \big( n \widehat{\mathcal{V}}_{12}^2  / 2 > x \big) \big\}^2, \\[.5em]
& b_2 = K(K-1)(K-2) \big\{ \mP \big( n \widehat{\mathcal{V}}_{12}^2  / 2 > x \big) \big\}^2 \quad \text{and} \quad b_3 = 0.
\end{align*}
Therefore it is enough to verify that $\mP ( n \widehat{\mathcal{V}}_{12}^2  / 2 > x ) = O(K^{-2})$ under the given conditions. Then we have $(I) \rightarrow 0$ as $n,K \rightarrow \infty$.

In what follows, we prove $\mP ( n \widehat{\mathcal{V}}_{12}^2  / 2 > x ) = O(K^{-2})$ and $(II) \rightarrow 0$. First we apply Theorem~\ref{Theorem: moderate deviation} with $x \asymp 4\lambda_1 \log K = o(n^\theta)$ to have
\begin{align*}
\frac{K(K-1)}{2} \mP \big( n \widehat{\mathcal{V}}_{12}^2  / 2 > x \big) = \frac{K(K-1)}{2} \mP \bigg( \sum_{v=1}^\infty \lambda_v  \xi_v^2  > x \bigg)\{1 + o(1)\}.
\end{align*}
Using the tail approximation given by \cite{zolotarev1961concerning} as $x \rightarrow \infty$:
\begin{align*}
&\mP \bigg( \sum_{v=1}^\infty \lambda_v  \xi_v^2  > x \bigg) = \frac{\kappa}{\Gamma(\mu_1/2)} \left( \frac{x}{2\lambda_1} \right)^{\mu_1/2-1} \exp\left( - \frac{x}{2\lambda_1} \right) \{1 + o(1)\},
\end{align*}
we have 
\begin{align*}
\frac{K(K-1)}{2} \mP \big( n \widehat{\mathcal{V}}_{12}^2  / 2 > x \big) & ~=~ \frac{\kappa}{\Gamma(\mu_1/2)} \left( \frac{x}{2\lambda_1} \right)^{\mu_1/2-1} \exp\left( - \frac{x}{2\lambda_1} \right) \{1 + o(1)\} \\[.5em]
&~=~   \exp\bigg\{ - \frac{2^{\mu_1/2-2}\kappa}{\Gamma(\mu_1/2)}  \exp\left( -\frac{y}{2}\right) \bigg\} \{1 + o(1)\}. 
\end{align*}
Therefore $\mP ( n \widehat{\mathcal{V}}_{12}^2  / 2 > x ) = O(K^{-2})$ and $(II) \rightarrow 0$ as $n,K \rightarrow \infty$, which completes the proof. 

\subsection{Proof of Theorem~\ref{Theorem: Upper Bound}}
Let us start by presenting some observations that are useful in the proof. 
\begin{itemize}
	\item \textbf{(O1).} From Lemma~\ref{Lemma: Unconditional Analysis}, we know that there exists a fixed constant $C_1 > 0$ such that 
	\begin{align*}
	\max_{1 \leq k < l \leq K}  \big| \widehat{\mathcal{V}}_{kl} - \mathcal{V}_{kl} \big| \leq C_1 \sqrt{\frac{B }{n_{\text{min}}} \log \left( \frac{K}{\beta} \right) }
	\end{align*}
	with probability at least $1-\beta$. \\[-0.5em]
 	\item \textbf{(O2).} Let us define $c_\alpha$ such that 
 	\begin{align} \label{definition: c_alpha}
 	c_\alpha := \inf \bigg\{ t \in \mathbb{R} : \frac{1}{N!}\sum_{\boldsymbol{b} \in \mathcal{B}_N} \mathds{1}\left( \widehat{\mathcal{V}}_{h,\text{max}}^{(\boldsymbol{b})} \geq t \right) \leq \alpha \bigg\}.
 	\end{align}
 	From Theorem~\ref{Theorem: K-sample concentration} and \textbf{(B2)}, there exists another fixed constant $C_2 >0$ such that 
 	\begin{align*}
 	c_\alpha \leq C_2 \sqrt{\frac{B }{n_{\text{min}}} \log \left( \frac{K}{\alpha} \right) }
 	\end{align*}
 	with probability one. Here we used the fact that $\widehat{\sigma}_K^2 \leq \max_{1 \leq i < j \leq N} \widetilde{h}(Z_i,Z_j) \leq B$. Thus the same inequality can be derived from (\ref{Eq: K-sample concentration}) in Theorem~\ref{Theorem: K-sample concentration} by replacing $\widehat{\sigma}_K^2$ with $\max_{1 \leq i < j \leq N} \widetilde{h}(Z_i,Z_j)$, which is more efficient to compute. \\[-0.5em]
 	\item \textbf{(O3).} Based on the definition of $c_\alpha$ in (\ref{definition: c_alpha}), observe that the event $\{p_{\text{perm}} > \alpha\}$, which is equivalent to
 	\begin{align*}
 	\bigg\{ \frac{1}{N!}\sum_{\boldsymbol{b} \in \mathcal{B}_N} \mathds{1}\left( \widehat{\mathcal{V}}_{h,\text{max}}^{(\boldsymbol{b})} \geq \widehat{\mathcal{V}}_{h,\text{max}} \right) > \alpha \bigg\},
 	\end{align*}
 	implies that $\{\widehat{\mathcal{V}}_{h,\text{max}} \leq c_\alpha\}$. 
\end{itemize}
Having these observations in mind, let us define an event $A_\beta$ such that
\begin{align*}
A_\beta = \Bigg\{ \max_{1 \leq k < l \leq K}  \big| \widehat{\mathcal{V}}_{kl} - \mathcal{V}_{kl} \big| \leq C_1 \sqrt{\frac{B }{n_{\text{min}}} \log \left( \frac{K}{\beta} \right) }  \Bigg\}.
\end{align*}
Then for sufficiently large $n_{\text{min}}$, the type II error of the permutation test is bounded by
\begin{align*}
& \mP \Big\{ p_{\text{perm}} > \alpha \Big\}  \\[.5em]
\overset{(i)}{\leq} ~ &  \mP \Big\{ \max_{1 \leq k < l \leq K}  \widehat{\mathcal{V}}_{kl} \leq c_\alpha \Big\}  \\[.5em]
\overset{(ii)}{\leq} ~ & \mP \Bigg\{ \max_{1 \leq k < l \leq K}  \widehat{\mathcal{V}}_{kl} \leq C_2 \sqrt{\frac{B}{n_{\text{min}}} \log \left( \frac{K}{\alpha} \right) } \Bigg\}  \\[.5em]
= ~ & \mP \Bigg\{ \max_{1 \leq k < l \leq K}  \widehat{\mathcal{V}}_{kl} \leq C_2 \sqrt{\frac{B}{n_{\text{min}}} \log \left( \frac{K}{\alpha} \right) } , ~ A_\beta \Bigg\}  + \mP \Bigg\{ \max_{1 \leq k < l \leq K}  \widehat{\mathcal{V}}_{kl} \leq C_2 \sqrt{\frac{B}{n_{\text{min}}} \log \left( \frac{K}{\alpha} \right) }, ~ A_\beta^c \Bigg\}  \\[.5em]
\overset{(iii)}{\leq} ~ & \mP \Bigg\{ \max_{1 \leq k < l \leq K}  \widehat{\mathcal{V}}_{kl} \leq C_2 \sqrt{\frac{B}{n_{\text{min}}} \log \left( \frac{K}{\alpha} \right) } , ~ A_\beta \Bigg\}  + \beta,
\end{align*}
where step~$(i)$ uses \textbf{(O3)}, step~$(ii)$ follows by \textbf{(O2)} and step~$(iii)$ uses \textbf{(O1)}. Furthermore, using the triangle inequality, we see that $\max_{1 \leq k < l \leq K}  \widehat{\mathcal{V}}_{kl} \geq \max_{1 \leq k < l \leq K}  \mathcal{V}_{kl} -  \max_{1 \leq k < l \leq K}  \big| \widehat{\mathcal{V}}_{kl} - \mathcal{V}_{kl} \big|$. Also note that $\max_{1 \leq k < l \leq K}  \mathcal{V}_{kl} \geq b_N r_N^\star$ under the given condition. Thus
\begin{align*}
&\mP \Bigg\{ \max_{1 \leq k < l \leq K}  \widehat{\mathcal{V}}_{kl} \leq C_2 \sqrt{\frac{B}{n_{\text{min}}} \log \left( \frac{K}{\alpha} \right) } , ~ A_\beta \Bigg\} \\[.5em]
\leq ~ & \mP \Bigg\{ b_N r_N^\star \leq C_1 \sqrt{\frac{B }{n_{\text{min}}} \log \left( \frac{K}{\beta} \right) } +  C_2 \sqrt{\frac{B}{n_{\text{min}}} \log \left( \frac{K}{\alpha} \right) } \Bigg\}.
\end{align*}
This gives an upper bound for the type II error that does not depend on $(P_1,\ldots,P_K)$. Since $B$ is constant under \textbf{(B1)}, the upper bound goes to zero by taking e.g. $\beta = 1/b_N$. This completes the proof.

\subsection{Proof of Corollary~\ref{Corollary: Upper Bound 2}}
First by the triangle inequality and Slutsky's argument,
\begin{align*}
\mP(p_{\text{MC}} > \alpha) & ~ \leq ~ \mP( |p_{\text{MC}} - p_{\text{perm}}| + p_{\text{perm}}  > \alpha)  \\[.5em]
& ~ \leq ~ \mP( |p_{\text{MC}} - p_{\text{perm}}| > \alpha/2) + \mP( p_{\text{perm}}  > \alpha/2).
\end{align*}
Followed by Theorem~\ref{Theorem: Upper Bound}, we have that 
\begin{align*}
\limsup_{n_{\text{min}} \rightarrow \infty} \sup_{(P_1,\ldots,P_K) \in \mathcal{F}_h(b_N r_N^\star)} \mP( p_{\text{perm}}  > \alpha/2) = 0.
\end{align*}
Therefore it suffices to control the first term. Let us write
\begin{align*}
F(t) = \frac{1}{N!}\sum_{\boldsymbol{b} \in \mathcal{B}_N} \mathds{1}\left( \widehat{\mathcal{V}}_{h,\text{max}}^{(\boldsymbol{b})} \leq  t \right) \quad \text{and} \quad F_M(t) = \frac{1}{M} \sum_{i=1}^M \mathds{1}\left( \widehat{\mathcal{V}}_{h,\text{max}}^{(\boldsymbol{b}_i^\prime)} \leq t \right).
\end{align*}
Then it can be shown that
\begin{align*}
|p_{\text{perm}} - p_{\text{MC}}| ~ \leq ~ \sup_{t \in \mathbb{R}} \big| F(t) - F_M(t) \big| + \frac{2}{M+1}.
\end{align*}
Hence the first term is bounded by
\begin{align*}
\mP( |p_{\text{MC}} - p_{\text{perm}}| > \alpha/2)  ~ \leq ~ \mP\left( \sup_{t \in \mathbb{R}} \big| F(t) - F_M(t) \big| > \frac{\alpha}{4} \right) + \mP\left( \frac{2}{M+1} > \frac{\alpha}{4} \right).
\end{align*}
Notice that by the DKW inequality \citep[e.g.][]{massart1990tight}, 
\begin{align*}
\mP\left( \sup_{t \in \mathbb{R}} \big| F(t) - F_M(t) \big| > \frac{\alpha}{4} \right)  \leq 2 e^{-M\alpha^2/8}.
\end{align*}
Thus 
\begin{align*}
\lim_{M \rightarrow \infty} \limsup_{n_{\text{min}} \rightarrow \infty} \sup_{(P_1,\ldots,P_K) \in \mathcal{F}_h(b_N r_N^\star)} \mP( |p_{\text{MC}} - p_{\text{perm}}| > \alpha/2) = 0,
\end{align*}
which results in the conclusion.

\subsection{Proof of Theorem~\ref{Theorem: Minimax Lower Bound}}
Motivated by Theorem 1 of \cite{tolstikhin2017minimax}, we use discrete distributions to prove the result. Specifically, we choose two distinct points $z_1,z_2$ on $\mathbb{R}^d$ such that $\varphi(0) - \varphi(z_1-z_2) \geq \kappa_1$. Consider the discrete distribution $p_0$ supported on the two points $z_1,z_2$ with probability $p_0(z_1) = 1/2$ and $p_0(z_2) = 1/2$. Consider another discrete distribution $p_1$ on the same support such that $p_1(z_1) = 1/2 + \delta$ and $p_1(z_2) = 1/2 - \delta$ where $\delta= b r_N^\star / \sqrt{2\kappa_1}$ and $b$ will be specified later. Then based on the translation invariant property of $h$, the MMD between $p_0$ and $p_1$ is calculated as
\begin{align} \label{Eq: Lower Bound of MMD}
\mathcal{V}_h(p_0,p_1) = \delta\sqrt{2 \{\varphi(0) - \varphi(z_1-z_2)\}} \geq \delta \sqrt{2\kappa_1}.
\end{align}
See \cite{tolstikhin2017minimax} for details. 

Next let $k$ be a discrete random variable uniformly distributed on $\{1,\ldots,K\}$. Then we set $P_{1,k}=p_0 \mathds{1}(k \neq 1) + p_1 \mathds{1}(k = 1), \ldots, P_{K,k}=p_0\mathds{1}(k \neq K) + p_1 \mathds{1}(k = K)$. Under this setting, it can be seen that $(P_{1,k},\ldots,P_{K,k}) \in \mathcal{F}_h(b r_N^\star)$ using (\ref{Eq: Lower Bound of MMD}). 

For each $k \in \{1,\ldots,K \}$, let $q_k$ be the joint probability function of $X_{1,1},\ldots,X_{n_K,K}$ given by
\begin{align*}
q_k(x_{1,1},\ldots,x_{n_K,K}) = & \prod_{i=1}^{n_1} \{ p_0(x_{i,1})\mathds{1}(k \neq 1) + p_1(x_{i,1})\mathds{1}(k=1)\} \times \\[.5em] 
& \cdots \times \prod_{i=1}^{n_K} \{ p_0(x_{i,K})\mathds{1}(k \neq K) + p_1(x_{i,K})\mathds{1}(k=K)\}. 
\end{align*}
Then we consider a mixture distribution given by $\overline{q}_{H_1} = \frac{1}{K} \sum_{k=1}^K q_k$. Also denote
\begin{align*}
q_{H_0}(x_{1,1},\ldots,x_{n_K,K}) = \prod_{i=1}^{n_1} p_0(x_{i,1}) \times \cdots \times \prod_{i=1}^{n_K} p_0(x_{i,K}).
\end{align*}
Then the likelihood ratio between $\overline{q}_{H_1}$ and $q_{H_0}$ is
\begin{align*}
L_N = ~ &  \frac{\overline{q}_{H_1}(X_{1,1},\ldots,X_{n_K,K})}{q_{H_0}(X_{1,1},\ldots,X_{n_K,K})} =\frac{1}{K} \sum_{k=1}^K  \prod_{i=1}^{n_k}  \frac{p_1(X_{i,k})}{p_0(X_{i,k})} = \frac{1}{K} \sum_{k=1}^K  \prod_{i=1}^{n_k}  \frac{p_0(X_{i,k}) + \delta\gamma(X_{i,k})}{p_0(X_{i,k})} \\[.5em]
= ~ & \frac{1}{K} \sum_{k=1}^K  \prod_{i=1}^{n_k}  \{ 1 + 2\delta \gamma(X_{i,k}) \},
\end{align*}
where $\gamma(X_{i,k}) = \{ \mathds{1}(X_{i,k} = z_1) - \mathds{1}(X_{i,k} = z_2)\}$. Moreover, the expected value of $L_N^2$ under $H_0$ is 
\begin{align*}
\mE_0 ( L_N^2 ) & = \frac{1}{K^2} \sum_{k=1}^K \sum_{k^\prime=1}^K  \mE_0 \Bigg[ \prod_{i=1}^{n_k}  \{ 1 + 2\delta \gamma(X_{i,k}) \} \prod_{i=1}^{n_k^\prime} \{ 1 + 2\delta \gamma(X_{i,k^\prime}) \} \Bigg] \\[.5em]
& = \frac{1}{K^2} \sum_{k=1}^K \mE_0 \Bigg[ \prod_{i=1}^{n_k}  \{ 1 + 2\delta \gamma(X_{i,k}) \}^2 \Bigg] \\[.5em]
& ~~~~~~+ \frac{1}{K^2} \sum_{k \neq k^\prime}^K \mE_0 \Bigg[ \prod_{i=1}^{n_k}  \{ 1 + 2\delta \gamma(X_{i,k}) \} \prod_{i=1}^{n_k^\prime} \{ 1 + 2\delta \gamma(X_{i,k^\prime}) \} \Bigg] \\[.5em]
& = \frac{1}{K^2} \sum_{k=1}^K  \prod_{i=1}^{n_k}  \{ 1 + 4\delta^2 \} + \frac{1}{K^2} \sum_{k \neq k^\prime}^K \prod_{i=1}^{n_k} \prod_{i=1}^{n_k^\prime} \{1\} \\[.5em]
& \leq \frac{1}{K^2} \sum_{k=1}^K \exp\left( 4 n_k \delta^2 \right) + \frac{K(K-1)}{K^2} 
\end{align*}
where the last inequality uses $1+ x \leq e^x$ for all $x$. From the assumption $\textbf{(B2)}$, we know that there exists a fixed constant $C_3 > 0$ such that 
\begin{align*}
\frac{1}{K^2} \sum_{k=1}^K \exp\left( 4 n_k \delta^2 \right)  \leq \frac{1}{K}  \exp\left( C_3 n_\text{min} \delta^2 \right).
\end{align*}
Finally, based on the standard $\chi^2$ method for minimax testing~\citep[e.g.][]{baraud2002non}, it is enough to find a positive constant $b$ such that $\delta = b r_N^\star / \sqrt{2\kappa_1} < 1/2$ and $\mE_0[L_N^2] \leq 1 +4(1-\alpha-\zeta)^2$. Indeed, this holds for any $b < \min \{ \sqrt{2\kappa_1/C_3}, \sqrt{\kappa_1} /(\sqrt{2}\kappa_2)\}$ for sufficiently large $K$, which completes the proof.

\end{document}